\documentclass[10pt]{article}

\usepackage{dirtytalk}

\usepackage{amsfonts}
\usepackage{amssymb}
\usepackage{amsmath}
\usepackage{enumerate}
\usepackage{amsthm}

\usepackage{mathtools}

\newcommand{\cf}{\mathrm{cf}}

\newcommand{\cov}{\mathrm{cov}}
\newcommand{\ran}{\mathrm{ran}}

\newcommand{\pp}{\mathrm{pp}}
\newcommand{\tcf}{\mathrm{tcf}}

\newtheorem{Th}{\bf THEOREM}[section]

\newtheorem{Pro}[Th]{\bf PROPOSITION}

\newtheorem{fact}[Th]{\bf FACT}

\newtheorem{Cor}[Th]{\bf COROLLARY}

\newtheorem{Obs}[Th]{\bf OBSERVATION}

\theoremstyle{definition} \newtheorem{Def}[Th]{\bf DEFINITION}
 
\theoremstyle{remark}

\usepackage[english]{babel}
\makeindex

\title{PIECE SELECTION AND CARDINAL ARITHMETIC}
\author{Pierre MATET}

\date{}

\begin{document}

\maketitle

\renewcommand{\thefootnote}{\arabic{footnote}} 	%Num屍otation avec les chiffres arabes exemple 1, 2, 3...

\renewcommand{\thefootnote}{}                                %  Aucun num屍o

%\begin{center}
%Department de Mathematics, University de Caen, 14032 CAEN cedex, France

%e-mail : matet@math.unicaen.fr
%\end{center}

\vskip 0,7cm

\begin{abstract}  
%The piece selection principle is one of those large cardinal like properties that can be satisfied at successor cardinals. 

We study the effects of piece selection principles on cardinal arithmetic (Shelah style). As an application, we discuss questions of Abe and Usuba. In particular, we show that if $\lambda \geq 2^\kappa$, then (a) $I_{\kappa, \lambda}$ is not $(\lambda, 2)$-distributive, and (b) $I_{\kappa, \lambda}^+ \rightarrow (I_{\kappa, \lambda}^+)^2_\omega$ does not hold.

%Let $\mu, \nu$ and $\kappa$ be three infinite cardinals such that $\cf(\mu) = \mu \in \nu \setminus \{\cf(\nu)\}$ and $\kappa = \nu^+ = 2^\nu$. Then by a result of Shelah \cite{She10}, $\diamondsuit_\kappa[J]$ holds for any normal ideal $J$ on $\kappa$ with $E^\kappa_\mu \in J^\ast$. By an earlier result of Gregory \cite{Gre76} and Shelah \cite{She78}, if we make the extra assumption that $\tau^\mu \leq \nu$ for every cardinal $\tau < \nu$, then even $\diamondsuit^\ast_\kappa (E^\kappa_\mu)$ will hold. We generalize these two results to $P_\kappa (\lambda)$.  
\end{abstract}

\bigskip

\section{Introduction}

\smallskip

In \cite{Abramson} Abramson, Harrington, Kleinberg and Zwicker pointed out that many large cardinal properties can be reformulated as flipping properties, which are of the following type : One is given a family $F$ of subsets of a set $X$. The property asserts that for some \say{flip} $h \in \prod_{A \in F} \{A, X \setminus A\}$, the $h (A)$'s satisfy some intersection property (for instance the finite intersection property). \say{Intersection} is taken in a wide sense so that e.g. \emph{diagonal} intersections are allowed. Notice that with the family $F$ is associated the family of all \emph{partitions} of the form $\{A, X \setminus A\} \setminus \{\emptyset\}$ for $A \in F$. A flip of $F$ can now be seen as a piece selection operation. Namely, for each partition $\{A, X \setminus A\}\setminus \{\emptyset\}$, we choose one piece, either $A$ or its complement. 

For a typical example, let $\kappa$ be a measurable cardinal, and $U$ be a ($\kappa$-complete) measure on $\kappa$. For $F$ take the collection of all partitions of $\kappa$ into one or two pieces. For each such partition, select the piece in $U$. Then the pieces chosen have the property that the intersection of any less than $\kappa$ many of them is cofinal in $\kappa$. Notice that since $U$ is $\kappa$-complete, it does not matter whether our partitions have one, two or more pieces, as long as the number of pieces is less than $\kappa$. By thus increasing the number of pieces, we obtain a natural generalization of the original flipping properties. In this extended framework, regularity of an infinite cardinal $\kappa$ can be expressed as the property that for any partition of $\kappa$ into less than $\kappa$ many pieces, one of the pieces must be cofinal in $\kappa$. The setting can be generalized by allowing $J$-partitions, and not just partitions. (Recall that for an ideal $J$ on a set $X$, a $J$-partition of $X$ is a subset $Q$ of $J^+$ such that
\begin{itemize}
\item $A \cap B \in J$ for any two distinct members $A, B$ of $Q$.
\item For any $C \in J^+$, there is $A \in Q$ with $A \cap C \in J^+$.)  
\end{itemize}
This is a way to handle properties defined in terms of distributivity. For a further generalization we relax the requirement that a piece has to be selected in \emph{each} partition. So for instance we're given $\kappa$ many partitions of $\kappa$, and we might be happy to pick one piece in $\kappa$ many partitions. Piece selection principles of this type have been introduced in our joint paper \cite{Laura} with Laura Fontanella. Their study is continued in the present paper. 

Our starting point is Solovay's celebrated result \cite{SOL} on strongly compact cardinals and the Singular Cardinal Hypothesis. This theorem can be revisited in a number of ways. For instance it is shown in \cite{pcf} that if $\cf (\lambda) < \kappa$ and there is a  weakly $\lambda^{++}$-saturated, $(\cf (\lambda))^+$-complete ideal on $P_\kappa (\lambda)$, then $\pp (\lambda) = \lambda^+$. In another direction, Usuba \cite{Usuba} established that if $\kappa$ is mildly $\lambda$-ineffable and $\cf (\lambda) \geq \kappa$, then $\lambda^{< \kappa} = \lambda$, or equivalently (since $\kappa$ is inaccessible), there is a cofinal subset of $P_\kappa (\lambda)$ of size $\lambda$ (i.e. $u (\kappa, \lambda) = \lambda$). Now it was noted \cite{Carlo}  from the very beginning that mild ineffability can be reformulated as a piece selection principle. We consider various weakenings of this principle and attempt to compare their relative strengths. Some of these properties can be satisfied at a weakly, but not strongly, inaccessible cardinal, or even at a successor cardinal. So this part of the paper is a contribution to the age-old program of determining what's left of weak or strong compactness when inaccessibility is removed. As the program developed, an impressive list of properties emerged, especially in connection with \emph{weak} compactness. Some of these properties can be tentatively classified as weak (the tree property), of medium strength (our $PS^+$) or strong (the weak compactness of the infinitary language $L_{\kappa\omega}$). For others (e.g. our $PS$), the situation is not so clear and further work is needed. Our version of Solovay's result reads as follows (see Proposition 3.4 and Observation 4.6).

\medskip

\begin{Th}  Suppose that $\cf (\lambda) < \kappa$ and $PS^+ ((\cf (\lambda))^+, \kappa, \lambda)$ holds. Then  $\cov (\lambda, \lambda, (\cf (\lambda))^+, 2) = \lambda^+$.
 \end{Th}
 
 \medskip

In the remainder of the paper, which is devoted to applications, we still deal with variants of mild ineffability, but this time inaccessibility of $\kappa$ is implied. It is a central problem in the theory of $P_\kappa (\lambda)$ to determine how the infinite Ramsey theorem generalizes in this framework. (Note that the theorem comes in several versions, from the weak \say{$\{\omega \} \rightarrow (I_\omega^+)^2_2$} to the strong \say{$I^+_\omega \rightarrow (I^+_\omega)^{n}_m$ for all finite $n, m$}). The study of partition relations on $P_\kappa (\lambda)$ is known to be tricky business. Carr \cite{Carr2} mentions that \say{repeated efforts to obtain} 
\begin{itemize}
\item $\{P_\kappa (\lambda) \} \rightarrow (I^+_{\kappa, \lambda})^2_2$ implies $\kappa$ is mildly $\lambda$-ineffable,  
\item $\kappa$ is mildly $\lambda$-ineffable implies $\{P_\kappa (\lambda) \} \rightarrow (I^+_{\kappa, \lambda})^3_2$
\end{itemize} 
\say{failed miserably}. Johnson asked in \cite{Johnson2} whether the $(\lambda, 2)$-distributivity of $I_{\kappa, \lambda}$ follows from the mild $\lambda$-ineffabilty of $\kappa$. This was answered in the negative by Abe \cite{Abe} who showed that if $\lambda^{< \kappa} = 2^\lambda$, then (a) $I_{\kappa, \lambda} \vert A$ is not $(\lambda, 2)$-distributive for any stationary $A$, and (b) $I^+_{\kappa, \lambda} \rightarrow (I^+_{\kappa, \lambda})^2_2$ does not hold. This led him to ask whether  $\lambda > \kappa$ implies that (a) $I_{\kappa, \lambda}$ is not $(\lambda, 2)$-distributive, and (b) $I^+_{\kappa, \lambda} \rightarrow (I^+_{\kappa, \lambda})^2_2$ fails. This was answered, again in the negative, by Shioya \cite{Shioya}. It should be noted that his model is obtained by adding many  Cohen subsets of $\kappa$. In fact it was shown in \cite{Cov} that if $\kappa$ is mildly $\lambda^{< \kappa}$-ineffable and $\lambda^{< \kappa} < \rm \bf{cov} (\rm \bf{M}_{\kappa, \lambda})$, then $I^+_{\kappa, \lambda} \rightarrow (I^+_{\kappa, \lambda})^n_\eta$ holds for any $n < \omega$ and any $\eta < \kappa$.  Since $\rm \bf{cov} (\rm \bf{M}_{\kappa, \lambda}) \leq \rm \bf{cov} (\rm \bf{M}_{\kappa, \kappa}) \leq \frak{d}_\kappa \leq 2^\kappa$, it made us think that maybe it could be proved that if $\lambda \geq 2^\kappa$, then (a) $I_{\kappa, \lambda}$ is not $(\lambda, 2)$-distributive, and (b) $I^+_{\kappa, \lambda} \rightarrow (I^+_{\kappa, \lambda})^2_2$ fails. 

\medskip

Part of the difficulty with the ordinary partition relation on $P_\kappa (\lambda)$ stems from the fact that $\subset$ is not a \emph{linear} ordering. To avoid this kind of pitfalls we chose to work with weaker partition relations (note that by negating them, we will obtain stronger results). Given an ideal $J$ on $P_\kappa (\lambda)$ and a coloring of $P_\kappa (\lambda) \times P_\kappa (\lambda)$, we are looking for a color $i$ and a set $A$ in $J^+$ that is not necessarily $i$-homogeneous (all pairs $(a, b)$ from $A \times A$ with $a \subset b$ have color $i$), but at least $i$-homogeneous mod $J$ (meaning that for each $a$ in $A$, the set of all $b$ in $A$ such that $(a, b)$ does not have color $i$ lies in $J$). We denote this particular partition relation by $\{P_\kappa (\lambda\} \xrightarrow{J} (J^+)^2_\rho$, where $\rho$ denotes the number of available colors. Still weaker partition relations are obtained in a similar fashion by coloring $\kappa \times P_\kappa (\lambda)$ (given $(a, b)$, we look at the color of $(\sup (a \cap \kappa), b)$) or $\kappa^+  \times P_\kappa (\lambda)$ (replace $\sup (a \cap \kappa)$ with $\sup (a \cap \kappa^+)$). To denote the corresponding partition property, we use $\xrightarrow[\kappa]{J}$ (respectively, $\xrightarrow[\kappa^+]{J})$. Proofs involve the usual ingredients ($\kappa$-normality, covering numbers, etc., and indeed some proofs are slight modifications of proofs of earlier results. Progress is achieved via a broader appeal to Shelah's pcf theory.

\medskip

Our efforts to prove the conjectures described above were only partially successful. By combining Observation 8.1 and Propositions 6.19, 8.4 and 9.6, one obtains the following.

\medskip

\begin{Th} Suppose that $2^\kappa \leq \lambda$, and let $D \in NS^\ast_{\kappa, \lambda}$. Then setting  $J = I_{\kappa, \lambda} \vert D$,  the following hold : 
 \begin{enumerate} [\rm (i)]
 \item  $J^+ \xrightarrow{J} (J^+)^2_\omega$ does not hold.
\item  $J^+ \xrightarrow{J} (J^+)^3_2$ does not hold.
 \item $J$ is not $(\lambda, 2)$-distributive.
 \end{enumerate}
\end{Th}

\medskip

There is a wide wide gap between this and what we can establish (see Proposition 5.32, Corollary 6.10, Observation 8.3 and Fact 9.5) under extra cardinal arithmetic assumptions such as Shelah's Strong Hypothesis (SSH).

\medskip

\begin{Th} Assuming SSH, the following hold :
 \begin{enumerate} [\rm (i)]
 \item  If $\overline{\mathfrak{d}}_\kappa \leq \lambda$ and $\cf (\lambda) \not= \kappa$, then for any $D \in NS^\ast_{\kappa, \lambda}$,  $I_{\kappa, \lambda} \vert D$ is not $(\kappa, 2)$-distributive and $(I_{\kappa, \lambda} \vert D)^+ \xrightarrow[\kappa]{I_{\kappa, \lambda}} (I_{\kappa, \lambda}^+, \omega_1)^2$ fails. If moreover $\kappa$ is weakly Mahlo, then for any $D \in NS^\ast_{\kappa, \lambda}$, $(I_{\kappa, \lambda} \vert D)^+ \xrightarrow[\kappa]{I_{\kappa, \lambda}} [I_{\kappa, \lambda}^+]^2_\lambda$ fails.
 \item  If $2^\kappa \leq \lambda$ and $\cf (\lambda) = \kappa$, then for any $D \in NS^\ast_{\kappa, \lambda}$,  $I_{\kappa, \lambda} \vert D$ is not $(\kappa^+, 2)$-distributive, and moreover $(I_{\kappa, \lambda} \vert D)^+ \xrightarrow[\kappa^+]{I_{\kappa, \lambda}} [I_{\kappa, \lambda}^+]^2_\lambda$ fails.
 \end{enumerate}
\end{Th} 

\medskip

On the positive side we have the following (see Corollary 9.2 and Fact 9.3).

\medskip

\begin{Th} Suppose that SSH holds and either $\cf (\lambda) = \kappa$, or $\cf (\lambda) < \kappa$ and $\lambda^+ < \mathfrak{d}_\kappa$, or $\cf (\lambda) > \kappa$ and $\lambda < \mathfrak{d}_\kappa$. Then the following hold : 
\begin{enumerate} [\rm (i)]
\item Suppose that $\kappa$ is weakly inaccessible. Then $I_{\kappa, \lambda}^+ \xrightarrow[\kappa]{I_{\kappa, \lambda}} [I_{\kappa, \lambda}^+]^2_{\kappa^+}$ holds. 
\item Suppose that $\kappa$ is weakly compact. Then $I_{\kappa, \lambda}$ is $(\kappa, 2)$-distributive, and moreover $I_{\kappa, \lambda}^+ \xrightarrow[\kappa]{I_{\kappa, \lambda}} (I_{\kappa, \lambda}^+)^2_\eta$ holds whenever $0 < \eta < \kappa$.
 \end{enumerate}
\end{Th} 

\medskip

Concerning $I_{\kappa, \lambda}^+ \rightarrow (I_{\kappa, \lambda}^+)^2_2$, it remains open whether it fails whenever $2^\kappa \leq \lambda$. What we do know is that it fails if $\lambda$ is large enough. In fact as shown in \cite{square}, $\{P_\kappa (\lambda)\} \xrightarrow{I_{\kappa, \lambda}} [I_{\kappa, \lambda}^+]^2_\lambda$ fails if $\lambda$ is large enough.

\medskip

The article is organized as follows. Section 2 is devoted to piece selection principles on $P_\kappa (\lambda)$, with emphasis on $PS^+ (\tau, \kappa, \lambda)$. It is shown that if $\cf (\lambda) < \kappa$ and $PS^+ ((\cf (\lambda))^+, \kappa, \lambda)$ holds, then there is no remarkably good scale on $\lambda$. Section 3 is concerned with piece selection principles on $\kappa$. It is observed that the tree property $TP (\kappa)$ is one of them. We use scales to establish that if $\cf (\lambda) < \kappa$, then $PS^+ ((\cf (\lambda))^+, \kappa, \lambda)$ implies $PS^+ ((\cf (\lambda))^+, \lambda^+)$. Section 4 is devoted to Shelah's covering numbers. It is shown that if $\lambda$ is singular and $PS^+ ((\cf (\lambda))^+, \lambda^+)$ holds, then $\cov (\lambda, \lambda, (\cf (\lambda))^+, 2) = \lambda^+$. In Section 5 we give cardinal arithmetic conditions under which for any club subset $C$ of $P_\kappa (\lambda)$, the partition property $(I_{\kappa, \lambda} \vert C)^+ \xrightarrow[\kappa]{ I_{\kappa, \lambda}} ( I_{\kappa, \lambda}^+, \rho)^2$ fails. This is continued in Section 6 where we deal with the stronger partition relations $(I_{\kappa, \lambda} \vert C)^+ \xrightarrow[\kappa]{ I_{\kappa, \lambda}} ( I_{\kappa, \lambda}^+)^2_2$ and $(I_{\kappa, \lambda} \vert C)^+ \xrightarrow{ I_{\kappa, \lambda}} ( I_{\kappa, \lambda}^+)^2_\omega$. Mild ineffability is the subject of Section 7. We prove that if $\kappa$ is mildly $\lambda$-ineffable and $\cf (\lambda) \not= \kappa$, then $\cov (\lambda, \kappa^+, \kappa^+,\kappa) = \lambda$. Section 8 contains results on the non-distributivity of $I_{\kappa, \lambda} \vert C$ for a club subset $C$ of $P_\kappa (\lambda)$. Finally in Section 9, we deal with the remaining case, that is the case when $\cf (\lambda) = \kappa$, and explain why this case must be handled separately.

\bigskip

\section{Piece selection}

\medskip

{\bf Throughout the paper $\kappa$ will denote a regular uncountable cardinal, and $\lambda$ a cardinal greater than or equal to $\kappa$}. We start with some definitions.

\medskip

\begin{Def}  For a set $A$ and a cardinal $\tau$,  we set $P_\tau (A) = \{ a \subseteq A :  \vert a \vert  < \tau\}$ and $[A]^\tau = \{x \subseteq A : \vert x \vert = \tau\}$.
\end{Def}

 \begin{Def} By a \emph{partition} of a set $X$ we mean a subset $Q$ of $P (X) \setminus \{\emptyset\}$ such that:
\begin{itemize}
\item $A \cap B = \emptyset$ for any two distinct members $A, B$ of $Q$.
\item $\bigcup Q = X$.
\end{itemize}
\end{Def}

 \begin{Def} An \emph{ideal} on a set $X$ is a nonempty collection $J$ of subsets of $X$ such that :
\begin{itemize}
\item $A \cup B \in J$ whenever $A,B \in J$.
\item $P(A) \subseteq J$ for all $A \in J$.
\item $X \notin J$.
%\item $\{x\}\in J$ for all $x\in X.$
\end{itemize}

Given an ideal $J$ on $X$, we denote by $J^+$ the set $\{A\subseteq X : A \notin J\}$, while $J^\ast$ denotes the set $\{A \subseteq X : X\setminus A \in J\}$. For any $A \in J^+$, we let $J \vert A = \{B \subseteq X : B \cap A \in J\}$.

We say that $J$ is \emph{$\kappa$-complete} if for any collection $Z$ of less than $\kappa$ many sets in $J$, one has $\bigcup Z \in J$. 
%The \emph{cofinality} of $J,$ denoted $Cof(J)$ is the least cardinality of a subcollection $I$ of $J$ such that $J=\bigcup_{A\in I} \mathcal{P}(A).$ For a cardinal $\tau,$ we say that $J$ is \emph{$\tau $-saturated} if there is no subset $Q\subseteq J^+$ of size $\tau$ with the property that $A\cap B\in J$  for every two distinct elements $A,B$ of $Q.$ $J$ is \emph{nowhere $\tau$-saturated} if there is no $A\in J^+$ such that $J\vert A$ is $\tau$-saturated. We say that $J$ is \emph{prime} if it is $2$-saturated. 

An ideal $K$ on $X$ \emph{extends} $J$ if $J \subseteq K$.

We let $I_\kappa = \bigcup_{\alpha < \kappa} P (\alpha)$ and
 
\centerline{$I_{\kappa,\lambda}= \bigcup_{a \in P_\kappa (\lambda)}P (\{b \in P_\kappa (\lambda) : a \setminus b \not= \emptyset\})$.} 

An ideal $J$ on $\kappa$ (respectively, $P_\kappa (\lambda)$) is \emph{fine} if it extends $I_\kappa$ (respectively, $I_{\kappa, \lambda}$). 

We let $NS_\kappa$ (respectively, $NS_{\kappa, \lambda}$) denote the nonstationary ideal on $\kappa$ (respectively, $P_\kappa (\lambda))$.
\end{Def}

 \begin{Def}  Let $\tau$ be an infinite cardinal less than or equal to $\kappa$. The piece selection principle $PS^+ (\tau, \kappa, \lambda)$ means that given a partition $Q_a$ of $P_\kappa (\lambda)$ with $\vert Q_a \vert < \tau$ for each $a \in P_\kappa (\lambda)$, there is $B \in I^+_{\kappa, \lambda}$ and $h \in \prod_{a \in P_\kappa (\lambda)} Q_a$ such that for any $a, b \in B$, the set $\{c \in h (a) \cap h (b) : a \cup b \subseteq c\}$ is nonempty.
 
$PS^\ast (\tau, \kappa, \lambda)$) (respectively, $PS (\tau, \kappa, \lambda)$) means that given a partition $Q_a$ of $P_\kappa (\lambda)$ with $\vert Q_a \vert < \tau$ for each $a \in P_\kappa (\lambda)$, we may find $B \in I^+_{\kappa, \lambda}$ and $h \in \prod_{a \in P_\kappa (\lambda)} Q_a$ such that for any $a, b \in B$, there is $t$ in $B$ (respectively, in $P_\kappa (\lambda)$) such that $a \cup b \subseteq t$ and the two sets $\{c \in h(a) \cap h(t) : t \subseteq c\}$ and $\{d \in h(b) \cap h(t) :  t \subseteq d\}$ are nonempty.
\end{Def}

\begin{Obs} $PS^+ (\tau, \kappa, \lambda) \Rightarrow PS^\ast (\tau, \kappa, \lambda) \Rightarrow PS (\tau, \kappa, \lambda)$.
\end{Obs}

\begin{Obs} Suppose that $PS^\ast (\kappa, \kappa, \lambda)$ holds. Then $\kappa$ is weakly inaccessible.
\end{Obs}

{\bf Proof.}  Suppose otherwise, and let $\kappa = \nu^+$. For $\nu \leq \gamma < \kappa$, select a bijection $j_\gamma : \gamma \rightarrow \nu$. Put $A = \{a \in P_\kappa (\lambda) : \nu \subseteq a\}$. For $a \in A$ and $i < \nu$, let $Q^{i}_a$ denote the collection of all $c \in P_\kappa (\lambda)$ such that $a \subseteq c$ and $j_{(\sup (c \cap \kappa)) + 1} (\sup (a \cap \kappa)) = i$. We may find $B \in I^+_{\kappa, \lambda} \cap P (A)$ and $i < \nu$ such that for any $a, b \in B$, there is $t$ in $B$ such that $a \cup b \subseteq t$ and the two sets $\{c \in Q^{i}_a \cap Q^{i}_t : t \subseteq c\}$ and $\{d \in Q^{i}_b \cap Q^{i}_t :  t \subseteq d\}$ are nonempty. Now pick $a, b \in B$ with $\sup (a \cap \kappa) < \sup (b \cap \kappa)$. There must be $t$ in $B$ with $a \cup b \subseteq t$, $c \in Q^{i}_a \cap Q^{i}_t$ with $t \subseteq c$ and $d \in Q^{i}_b \cap Q^{i}_t$ with  $t \subseteq d$. But then 
\begin{itemize}
\item $j_{(\sup (c \cap \kappa)) + 1} (\sup (a \cap \kappa)) = j_{(\sup (c \cap \kappa)) + 1} (\sup (t \cap \kappa))$.
\item $j_{(\sup (d \cap \kappa)) + 1} (\sup (b \cap \kappa)) = j_{(\sup (d \cap \kappa)) + 1} (\sup (t \cap \kappa))$.
\end{itemize}
It follows that $\sup (a \cap \kappa) = \sup (t \cap \kappa) = \sup (b \cap \kappa)$. Contradiction.
\hfill$\square$

%\begin{fact}  $IE_\kappa^2(\mathcal{I}_{\kappa,\lambda})$ implies $PS^+(\kappa,\lambda).$ \end{fact}
 
\begin{Obs} Suppose that $PS^+ (\tau, \kappa, \lambda)$ holds. Let $A \in I^+_{\kappa, \lambda}$, and for each $a \in A$, let $Q_a$ be a partition of the set $\{c \in A : a \subseteq c\}$ with $\vert Q_a \vert < \tau$. Then there is $B \in I^+_{\kappa, \lambda} \cap P (A)$ and $h \in \prod_{a \in B} Q_a$ such that for any $a, b \in B$, $h (a) \cap h (b) \not= \emptyset$. 
\end{Obs}

{\bf Proof.}  Define $\psi : P_\kappa (\lambda) \rightarrow A$ so that $x \subseteq \psi (x)$ for all $x \in P_\kappa (\lambda)$. For $x \in P_\kappa (\lambda)$, let $T_x$ denote the set of all $z \in P_\kappa (\lambda)$ such that $x \subseteq z$ but $\psi (x) \setminus z \not = \emptyset$, and set

\centerline{$Z_x = \{\{z \in P_\kappa (\lambda) : \psi (x) \subseteq z$ and $\psi (z) \in W \} : W \in Q_{\psi (x)}\}$.}

We may find $H \in I^+_{\kappa, \lambda}$ and $k \in \prod_{x \in P_\kappa (\lambda)} (Z_x \cup \{T_x\})$ such that $k (x) \cap k (y) \not= \emptyset $ for any $x, y \in H$. 

 \medskip

{\bf Claim.}  Let $x \in H$. Then $k (x) \in Z_x$.

\smallskip

{\bf Proof of the claim.} Suppose otherwise. Pick $y \in H$ with $\psi (x) \subseteq y$, and $z \in k (x) \cap k (y)$. Then $\psi (x) \subseteq y \subseteq z$. This contradiction completes the proof of the claim. 

\medskip

Now put $B = \psi``H$, and define $f : B \rightarrow H$ such that $\psi (f (a)) = a$ for all $a \in B$. Notice that $B \in I^+_{\kappa, \lambda} \cap P (A)$. Let $h \in \prod_{a \in B} Q_a$ be such that for any $a \in B$,

\centerline{$k (f (a)) = \{z \in P_\kappa (\lambda) : \psi (f (a)) \subseteq z$ and $\psi (z) \in h (a) \}$.}

Given $a, b \in B$, we may find $z$ in $k (f (a)) \cap k (f (b))$. Then $a \cup b = \psi (f (a)) \cup \psi (f (b)) \subseteq z \subseteq \psi (z)$, and moreover $\psi (z) \in h (a) \cap h (b)$. 
\hfill$\square$

\begin{Obs} Suppose that $PS^+ (\tau, \kappa, \lambda)$ holds. Then for any cardinal $\chi$ with $\kappa \leq \chi < \lambda$, $PS^+ (\tau, \kappa, \chi)$ holds. 
\end{Obs}

{\bf Proof.}  Let $\chi$ be a cardinal with $\kappa \leq \chi < \lambda$, and for each $y \in P_\kappa (\chi)$, let $Q_y$ be a partition of $P_\kappa (\chi)$ with $\vert Q_y \vert < \tau$. For $a \in P_\kappa (\lambda)$, put

\centerline{$W_a = \{\{c \in P_\kappa (\lambda) : c \cap \chi \in Z\} : Z \in Q_{a \cap \chi}\}$.}

Note that $\{c \cap \chi : c \in S\} \in Q_{a \cap \chi}$ for all $S \in W_a$. We may find $B \in I^+_{\kappa, \lambda}$ and $h \in \prod_{a \in P_\kappa (\lambda)} W_a$ such that for any $a, b \in B$, 
 
 \centerline{$\{c \in h (a) \cap h (b) : a \cup b \subseteq c\} \not= \emptyset$.}
 
 Set $Y = \{b \cap \chi : b \in B\}$. Notice that $Y \in I^+_{\kappa, \chi}$. Select $\psi : Y \rightarrow B$ so that for any $y \in Y$, $y = \psi (y) \cap \chi$, and define $k \in \prod_{y \in Y} Q_y$ by $k (y) = \{c \cap \chi : c \in h (\psi (y))\}$. Now given $x, y \in Y$, pick $c \in h (\psi (x)) \cap h (\psi (y))$ with  $\psi (x) \cup \psi (y) \subseteq c$. Then clearly, $x \cup y \subseteq c \cap \chi$, and moreover $c \cap \chi \in k (x) \cap k (y)$.   
 \hfill$\square$
 
 \medskip
 
 Let us next recall some material concerning scales in pcf theory.
 
 \medskip

\begin{Def} Let $A$ be an infinite set of regular cardinals such that $\vert A \vert < \min A$, and $I$ be an ideal on $A$ such that $\{A \cap a : a \in A \} \subseteq I$. 

We let $\prod A = \prod_{a \in A} a$. For $f, g \in \prod A$, we let $f <_I g$ if $\{a \in A : f (a) \geq g (a) \} \in I$. 

Let $\pi$ be a regular cardinal greater than $\sup A$. An increasing, cofinal sequence $\vec{f} = \langle f_\alpha : \alpha < \pi \rangle$  in $(\prod A, <_I)$ is said to be a \emph{scale of length} $\pi$. If there is such a sequence, we set $\tcf (\prod A / I) = \pi$. 
\end{Def}

\begin{fact} {\rm (\cite[Theorem 1.5 p. 50]{SheCA})}  Suppose that $\lambda$ is a singular cardinal. Then there is a set $A$ of regular cardinals such that $o.t. (A) = \cf (\lambda) < \min A, \sup A = \lambda$ and $\tcf (\prod A/I) =  \lambda^+$,  where $I$ is the noncofinal ideal on $A$. 
\end{fact}

\begin{Def} Let $\vec{f} = \langle f_\alpha : \alpha < \pi \rangle$ be an increasing, cofinal sequence in $(\prod A, <_I)$. An infinite limit ordinal $\delta < \pi$  is a \emph{good} (respectively, \emph{remarkably good}) \emph{point} for $\vec f$ if there is a cofinal (respectively, closed unbounded) subset $X \subseteq \delta$, and $Z_\xi \in I$ for $\xi \in X$ such that $f_\beta (a) < f_\xi (a)$ whenever $\beta < \xi$ are in $X$ and $a \in A \setminus (Z_\beta \cup Z_\xi)$.  
$\delta$ is a \emph{better point} for $\vec f$ if we may find a closed unbounded subset $X$ of $\delta$, and $Z_\xi \in I$ for $\xi\in X$ such that $f_\beta (a) < f_\xi (a)$ whenever $\beta < \xi$ are in $X$ and $a \in A \setminus Z_\xi$. 
$\delta$ is a \emph{very good point} for $\vec f$ if there is a closed unbounded subset $X$ of $\delta$, and $Z\in I$ such that $f_\beta (a) < f_\xi (a)$ whenever $\beta < \xi$ are in $X$ and $a \in A \setminus Z$. 

 The scale $\vec f = \langle f_\alpha : \alpha < \pi \rangle$  is \emph{good} (respectively, \emph{remarkably good}, \emph{better}, \emph{very good})  if there is a closed unbounded subset $C$ of $\pi$  with the property that every infinite limit ordinal $\delta$ in $C$ such that $\cf (\delta) < \sup A$ and $I$ is not $\cf (\delta)$-complete is a good (respectively, remarkably good, better, very good) point for $\vec f$.
\end{Def}

%\begin{Obs} \begin{enumerate}[\rm (i)]\item very good $\Rightarrow$ better $\Rightarrow$ remarkably good $\Rightarrow$ good.\item If $\delta$ is good, then $\cf(\delta) < \sup A$.\item If $\delta$ is good (respectively, remarkably good) and of uncountable cofinality, then there is a closed unbounded subset $C$ of $\delta$ consisting of good (respectively, remarkably good) points. \end{enumerate}\end{Obs} 

\begin{fact}  {\rm(\cite{CFM}, \cite{Norm})} Let $\delta < \pi$  be an infinite limit ordinal such that $I$ is $\cf (\delta)$-complete (respectively, $(\cf(\delta))^+$-complete). Then $\delta$  is a better (respectively, very good) point for $\vec f$.
\end{fact}

\medskip

We will show that if $\cf (\lambda) < \kappa$ and $PS^+ ((\cf (\lambda))^+, \kappa, \lambda^+)$ holds, then there is no remarkably good scale on $\lambda$.

\medskip

\begin{Def} Given two infinite cardinals $\tau$ and $\chi$ such that $\tau \leq \chi = \cf(\chi)$, we let $E^\chi_\tau$ (respectively, $E^\chi_{< \tau}$) denotes the set of all infinite limit ordinals $\alpha < \chi$ such that $\cf(\alpha) = \tau$ (respectively, $\cf(\alpha) < \tau$).
\end{Def}

\begin{Obs}  Suppose that $\tcf (\prod A  / I) = \pi$,  where
\begin{itemize}
\item $A$ is an infinite set of regular cardinals such that $\vert A \vert < \min A$,
\item $I$  is an ideal on $A$ such that $\{ A \cap a : a \in A\} \subseteq I,$
\end{itemize}
and $\vec f = \langle f_\alpha : \alpha < \pi\rangle$ is an increasing, cofinal sequence in $(\prod A, <_I)$. Let $\chi$ be an infinite cardinal with $\chi \leq (\sup A)^+$. Then the following are equivalent :
\begin{enumerate} [\rm (i)]
\item There is  a closed unbounded subset  $C$ of $\pi$ such that for any regular infinite cardinal $\theta < \chi$,  and any $\delta\in C\cap E_\theta^\pi$,  $\delta$  is a remarkably good point for $\vec f$. 
\item There  is  a closed unbounded subset  $D$ of $\pi$ such that  for any $e \in P_{\chi} (D)$,  there is $g : e \rightarrow I$  such that  $f_\alpha (a) < f_\beta (a)$  whenever $\alpha < \beta$  are in $e$ and $a \in A \setminus ( g(\alpha) \cup g(\beta))$. 
\end{enumerate}
\end{Obs}

{\bf Proof.}  \hskip0,4cm  (i) $\rightarrow$ (ii) :  By Proposition 8.11 of \cite{Scales} (we can take $D = C$).

\medskip

\hskip0,2cm  (ii) $\rightarrow$ (i) :  Easy (take $C =$ the set of limit points of $D$).
\hfill$\square$

\begin{Pro} Suppose that $\tcf (\prod A  / I) = \lambda^+$,  where
\begin{itemize}
\item $A$ is a set of regular cardinals such that $\vert A \vert < \min \{\kappa, \min A\}$ and $\sup A = \lambda$.
\item $I$  is an ideal on $A$ such that $\{ A \cap a : a \in A\} \subseteq I.$
\item $PS^+ (\vert A \vert^+, \kappa, \lambda^+)$ holds.
\end{itemize}
Let $\vec f = \langle f_\alpha : \alpha < \lambda^+ \rangle$ be an increasing, cofinal sequence in $(\prod A, <_I)$, and let $S$ denote the set of all infinite limit points $\delta < \lambda^+$ such that
\begin{itemize}
\item $\cf (\delta) < \kappa$.
\item $\delta$ is not a remarkably good point for $\vec f$.
\end{itemize}
Then $S$ is stationary in $\lambda^+$.
\end{Pro}

{\bf Proof.} Suppose otherwise, and select a closed unbounded subset $C$ of $\lambda^+$ such that any infinite limit ordinal $\delta$ in $C$ of cofinality less than $\kappa$ is a remarkably good point for $\vec f$. By Fact 2.14, for any $e\in P_\kappa (C)$,  there is $g_e : e \rightarrow I$  such that  $f_\alpha (a) < f_\beta (a)$  whenever $\alpha < \beta$  are in $e$ and $a \in A \setminus ( g_e (\alpha) \cup g_e (\beta))$. Pick a bijection $k : \lambda^+ \rightarrow C$ and for each nonempty $b \in P_\kappa (\lambda^+)$, $t_b \in \prod_{\alpha \in k``b} (A \setminus g_{k``b} (\alpha))$. Put $X = \{x \in P_\kappa (\lambda^+) : \sup x \in x\}$. For $x \in X$ and $a \in A$, set 

\centerline{$Q^{a}_x = \{ b \in X : x \subseteq b$ and $t_b (k (\sup x)) = a\}.$}

By Observation 2.7 we may find $B \in I^+_{\kappa, \lambda} \cap P (X)$ and $h : B \rightarrow A$ such that for any $x, y \in B$, $Q^{h (x)}_x \cap Q^{h (y)}_y \not= \emptyset$. There must be $H \in I^+_{\kappa, \lambda^+} \cap P (B)$ and $a \in A$ such that $h$ takes the constant value $a$ on $H$. Put $D = \{k (\sup x) : x \in H\}$.   

\medskip

{\bf Claim.}  Let $\alpha < \beta$ in $D$. Then $f_\alpha (a) < f_\beta (a)$. 

\smallskip

{\bf Proof of the claim.} Pick $x, y \in H$ with $k (\sup x) = \alpha$ and $k (\sup y) = \beta$, and $b \in Q^a_x \cap Q^a_y$. Then $t_b (\alpha) = a = t_b (\beta)$. Thus $a \in A \setminus (g_{k``b} (\alpha) \cup g_{k``b} (\beta))$, and therefore $f_\alpha (a) < f_\beta (a)$, which completes the proof of the claim.

\medskip

By the claim, the function $u : D \rightarrow A$ defined by $u (\alpha) = f_\alpha (a)$ is one-to-one. Contradiction.
\hfill$\square$

\medskip

Let us observe that by replacing the hypothesis that $PS^+ (\vert A \vert^+, \kappa, \lambda^+)$ holds with the stronger hypothesis that $\kappa$ is mildly $\lambda^+$-ineffable, we can actually prove \cite{Norm} that for \emph{cofinally many} regular uncountable cardinals $\sigma < \kappa$, the set of all points $\delta \in E^{\lambda^+}_\sigma$ that are not remarkably good for $\vec f$ is stationary.   

The conclusion of Proposition 2.15 entails the failure of a square principle of the following type.

\medskip

\begin{Def} Let $\theta$ and $\chi$ be two infinite cardinals such that $\theta \leq \chi = \cf(\chi)$, and $S \subseteq E^\chi_{< \chi}$. Then ${\rm WWS}^\theta_\chi (S)$ asserts the existence of $\mathcal{C_\gamma}$ for $\gamma < \chi$ such that
\begin{itemize}
\item $\vert \mathcal{C_\gamma} \vert < \chi$ ;
\item $\mathcal{C_\gamma}\subseteq P_\theta (\gamma)$ ;
\item if $\alpha \in S$, then there is a closed unbounded subset $C$ of $\alpha$ of order type $\cf(\alpha)$ such that $C \cap \gamma \in \bigcup_{ D \in \mathcal{C_\gamma}} P(D)$ for every $\gamma \in C$.  
\end{itemize}
\end{Def}

\begin{fact}   {\rm(\cite{Scales})} Let $A$, $I$ and $\pi$ be such that
\begin{itemize}
\item $A$ is an infinite set of regular cardinals ;
\item $\vert A \vert < \min A $ ;
\item $\sup A < \pi$ ;
\item $I$ is an ideal on $A$ such that $\{A \cap a : a \in A \} \subseteq I$ ;
\item  $\tcf (\prod A / I) = \pi$.
\end{itemize}
Further let $T$ be a collection of regular cardinals such that for any $\sigma \in T$, 
\begin{itemize}
\item $\sigma < \sup A$ ;
\item ${\rm WWS}^{\sup A}_\pi (E^\pi_\sigma)$ holds.
\end{itemize}
Then there is an increasing, cofinal sequence $\vec{f} = \langle f_\alpha : \alpha < \pi \rangle$ in $(\prod A, <_I)$ such that for any $\sigma \in T$, and any $\zeta \in E^\pi_\sigma$, $\zeta$  is a better point for $\vec{f}$.
\end{fact}

\begin{Pro} Suppose that $\cf (\lambda) < \kappa$ and $PS^+ ((\cf (\lambda))^+, \kappa, \lambda^+)$ holds. Then ${\rm WWS}^{\lambda}_{\lambda^+} (E^{\lambda^+}_\sigma)$ fails for some regular cardinal $\sigma$ with $\cf (\lambda) < \sigma < \kappa$.
 \end{Pro}

{\bf Proof.} By Proposition 2.15 and Facts 2.10, 2.12 and 2.17.
\hfill$\square$

\medskip

We next consider the two-cardinal version of the tree property.

\medskip

%\begin{Def}  A subset $C$ of $\mathcal{P}_\kappa(\lambda)$ is \emph{strongly closed} if $\bigcup X \in C$ for all $X \subseteq C$ with $0 < \vert X \vert < \kappa$.We let $SNS_{\kappa, \lambda}$ denote the collection of all $B \subseteq P_\kappa (\lambda)$ such that $B \cap C = \emptyset$ for some strongly closed, cofinal subset of $P_\kappa (\lambda)$. \end{Def}\begin{fact} {\rm((\cite{Men74}, \cite{Matet89})\begin{enumerate} [\rm (i)]\item $SNS_{\kappa,\lambda}$ is a $\kappa$-complete fine ideal on $P_\kappa (\lambda)$.\item $SNS_{\kappa,\lambda}$ is the collection of all $B \subseteq P_\kappa (\lambda)$ such that $\{a \in P_\kappa (\lambda) : f`` a \subseteq P (a)\}\cap B =\emptyset$ for some $f : \lambda \rightarrow P_\kappa (\lambda)$.\item Let $j$ be a one-to-one function from $\lambda \times \lambda$ to $\lambda$, and $D$ be the set of all $a \in P_\kappa (\lambda)$ such that $j (\alpha, \beta) \in a$ whenever $\alpha < \beta$ are in $a$. Then $NS_{\kappa, \lambda} = SNS_{\kappa, \lambda} \vert D$.\end{enumerate}\end{fact}

\begin{Def} $TP (\kappa,\lambda)$ asserts the following. Let $s_a \subseteq a$ for $a \in P_\kappa (\lambda)$ be such that for some $C \in NS_{\kappa,\lambda}^*$, we have $\vert \{s_a \cap c : c \subseteq a\} \vert <\kappa$ for all $c \in C$. Then there is $S \subseteq \lambda$ with the property that for every $b \in P_\kappa (\lambda)$, there is $a \in P_\kappa (\lambda)$ such that $b \subseteq a$ and $S \cap b= s_a \cap b$.
\end{Def}

\begin{fact} {\rm((\cite{Weiss})} It is consistent relative to a supercompact cardinal that $TP(\omega_2, \chi)$ holds for every cardinal $\chi \geq \omega_2$. 
\end{fact}

\begin{Obs} Let $s_a \subseteq a$ for $a \in P_\kappa (\lambda)$. Suppose that there is $C \in I_{\kappa,\lambda}^+$ such that $\vert \{s_a \cap c : c \subseteq a\} \vert <\kappa$ for all $c \in C$. Then $\vert \{s_a \cap d : d \subseteq a\} \vert <\kappa$ for all $d \in P_\kappa (\lambda)$. 
\end{Obs} 

{\bf Proof.} Suppose otherwise. Then we may find $d \in P_\kappa (\lambda)$ and $a_i \in P_\kappa (\lambda)$ for $i < \kappa$ so that
\begin{itemize}
\item $d \subseteq a_i$ for all $i < \kappa$.
\item $s_{a_i} \cap d \not= s_{a_j} \cap d$ whenever $i < j < \kappa$.
\end{itemize}
Pick $c \in C$ with $d \subseteq c$, and let $i < j < \kappa$. Then $s_{a_i} \cap c \not= s_{a_j} \cap c$, since otherwise we would have $s_{a_i} \cap d = s_{a_j} \cap d$. Contradiction. 
\hfill$\square$ 

\begin{fact} $PS(\kappa, \kappa,\lambda)$ implies $TP (\kappa,\lambda)$.
\end{fact}

{\bf Proof.} By Proposition 5.4 of \cite{Laura} and Observation 2.20.
\hfill$\square$ 

\medskip

We will now see that if in the definition of $PS^+$ we insist on selecting a piece from \emph{every} partition, then what we obtain is an apparently much stronger principle.

\medskip

 \begin{Def}  For a fine ideal $J$ on $P_\kappa (\lambda)$, the \emph{ideal extension principle} $IE (\tau, \kappa, \lambda, J)$ means that given a partition $Q_a$ of $P_\kappa (\lambda)$ with $\vert Q_a \vert < \tau$ for each $a \in P_\kappa (\lambda)$, there is $h \in \prod_{a \in P_\kappa (\lambda)} Q_a$ and an ideal $K$ on $P_\kappa (\lambda)$ extending $J$ such that $ran (h) \subseteq K^\ast$.
\end{Def}

\begin{Obs} $IE (\omega, \kappa, \lambda, J)$ holds.
\end{Obs} 

{\bf Proof.} Let $Q_a$ be a finite partition of $P_\kappa (\lambda)$ for $a \in P_\kappa (\lambda)$. Select a prime ideal $K$ extending $J$. For each $a \in P_\kappa (\lambda)$, there must be $W_a \in Q_a$ with $W_a \in K^\ast$. Now define $h \in \prod_{a \in P_\kappa (\lambda)} Q_a$ by $h (a) = W_a$.
\hfill$\square$ 

\begin{fact} {\rm(\cite{MPS1}} There is a partition of $P_\kappa (\lambda)$ into $\lambda^{< \kappa}$ sets in $I^+_{\kappa, \lambda}$.
 \end{fact}
 
%{\bf Proof.} Let $\langle a_\alpha : \alpha < \lambda^{< \kappa} \rangle$ be a one-to-one enumeration of $P_\kappa (\lambda)$. Select a bijection $j : \lambda^{< \kappa} \times \lambda^{< \kappa} \rightarrow \lambda^{< \kappa}$. Inductively define a one-to-one sequence $\langle b_\gamma : \gamma < \lambda^{< \kappa} \rangle$ of members of $P_\kappa (\lambda)$ so that if $\gamma = j (\alpha, \xi)$, then $a_\xi \subseteq b_\gamma$. Now for each $\alpha < \lambda^{< \kappa}$, set $A_\alpha = \{ b_{j (\alpha, \xi)} : \xi < \lambda^{< \kappa}\}$. \hfill$\square$ 

\begin{Obs} The following are equivalent :
\begin{enumerate} [\rm (i)]
\item $IE (\tau, \kappa, \lambda, I_{\kappa, \lambda})$ holds. 
\item Given a partition $Q_a$ of $P_\kappa (\lambda)$ with $\vert Q_a \vert < \tau$ for each $a \in P_\kappa (\lambda)$, there is $h \in \prod_{a \in P_\kappa (\lambda)} Q_a$ such that for any $a, b \in P_\kappa (\lambda)$,  there is $c \in h(a) \cap h(b)$ with $a \cup b \subseteq c$.
% : t \subseteq c\}$ such that $a\cup b \subseteq t$ and the two sets $\{c \in h(a) \cap h(t) : t \subseteq c\}$ and $\{d \in h(b) \cap h(t) :  t \subseteq d\}$ are nonempty.
\end{enumerate}
\end{Obs} 

{\bf Proof.}  \hskip0,4cm  (i) $\rightarrow$ (ii) : Trivial.

\medskip

\hskip0,2cm  (ii) $\rightarrow$ (i) :  Assume that (ii) holds, and let $Q_a$ be a partition of $P_\kappa (\lambda)$ with $\vert Q_a \vert < \tau$ for each $a \in P_\kappa (\lambda)$. By Fact 2.25, we may find a partition $T$ of $P_\kappa (\lambda)$ into $\lambda^{< \kappa}$ sets in $I^+_{\kappa, \lambda}$. Let $\langle T_w : w \in P_\kappa (\lambda) \rangle$ be a one-to-one enumeration of $T$. Pick a bijection $F : P_\kappa (\lambda) \rightarrow P_\omega (P_\kappa (\lambda)) \setminus \{ \emptyset\}$. For $w \in P_\kappa (\lambda)$ and $x \in T_w$, put $W_x  = \{ \bigcap_{a \in F (w)} k (a) : k \in \prod_{a \in F (w)} Q_a \}$. We may find $g \in \prod_{x \in P_\kappa (\lambda)} W_x$ such that for any $x, y \in P_\kappa (\lambda)$, there is $z \in g(x) \cap g (y)$ with $x \cup y \subseteq z$.
%$t \in P_\kappa (\lambda)$ such that $x \cup y \subseteq t$ and the two sets $\{c \in g(x) \cap g (t) : t \subseteq c\}$ and $\{d \in g (y) \cap g (t) :  t \subseteq d\}$ are nonempty. 
For $w \in P_\kappa (\lambda)$ and $x \in T_w$, let $k_x \in \prod_{a \in F (w)} Q_a $ be such that $g (x) =  \bigcap_{a \in F (w)} k_x (a)$. 
 
 \medskip

{\bf Claim 1.}  Let $v, w \in P_\kappa (\lambda)$, $x \in T_v$, $y \in T_w$ and $a \in F (v) \cap F (w)$. Then $k_x (a) = k_y (a)$. 

\smallskip

{\bf Proof of Claim 1.} Suppose otherwise. Then $k_x (a) \cap k_y (a) = \emptyset$. Since $g (x) \subseteq k_x (a)$ and $g (y) \subseteq k_y (a)$, it follows that $g (x) \cap g (y) = \emptyset$. This contradiction completes the proof of the claim.

\medskip

Put $h = \bigcup_{x \in P_\kappa (\lambda)} k_x$. Using Claim 1, it is easy to see that $h \in \prod Q_a$.

 \medskip

{\bf Claim 2.}  Let $e \in P_\omega (P_\kappa (\lambda)) \setminus \{ \emptyset\}$. Then $\bigcap_{a \in e} h (a) \in I^+_{\kappa, \lambda}$. 

\smallskip

{\bf Proof of Claim 2.} Let $e = F (w)$. Now given $s \in P_\kappa (\lambda)$, pick $x \in T_w$ with $s \subseteq x$. There must be $z \in \bigcap_{a \in F (w)} k_x (a)$ with $x \subseteq z$. Then clearly, $s \subseteq z$, and moreover $z \in \bigcap_{a \in e} h (a)$. This completes the proof of the claim and that of the observation.
 \hfill$\square$

\bigskip

\section{Piece selection at $\kappa$}

\medskip

In this section we concentrate on the case $\lambda = \kappa$.

\begin{Def} Given an infinite cardinal $\tau$, we let $PS^+(\tau, \kappa)$ assert the following: For $\beta \in \kappa$, let $Q_\beta$ be a partition of $\kappa \setminus \beta$ with $\vert Q_\beta \vert <\tau$. Then there is a cofinal subset $B$ of $\kappa$ and $h \in \prod_{\beta\in B} Q_\beta$ such that for any $\alpha, \beta \in B$, we have $h(\alpha) \cap h(\beta) \not= \emptyset$.

$PS^\ast(\tau, \kappa)$ (respectively $PS(\tau, \kappa)$) asserts the following: For $\beta \in \kappa$, let $Q_\beta$ be a partition of $\kappa \setminus \beta$ with $\vert Q_\beta \vert <\tau$. Then we may find a cofinal subset $B$ of $\kappa$ and $h \in \prod_{\beta \in \kappa} Q_\beta$ such that for any $\alpha, \beta \in B$, there is $\zeta$ in $B$ (respectively, in $\kappa$) such that $\max \{\alpha,\beta\} \leq \zeta$ and we have $h(\alpha) \cap h(\zeta) \not= \emptyset$ and $h(\beta)\cap h(\zeta) \not= \emptyset$.
%$PS^+(\kappa, \kappa)$ (respectively $PS^\ast(\kappa, \kappa)$, $PS(\kappa, \kappa)$) is abbreviated as $PS^+(\kappa)$ (respectively $PS^\ast(\kappa)$, $PS(\kappa)$) . 
\end{Def}

\begin{Obs} The following are equivalent :
 \begin{enumerate} [\rm (i)]
 \item $PS^+ (\tau, \kappa)$.
 \item $PS^+ (\tau, \kappa, \kappa)$.
 \end{enumerate}
\end{Obs} 

{\bf Proof.}  \hskip0,4cm  (i) $\rightarrow$ (ii) : Suppose that (i) holds. For $b \in P_\kappa (\kappa)$, let $Q_b$ be a partition of the set $\{c \in P_\kappa (\kappa) : b \subseteq c\}$ into less than $\tau$ many pieces. For $b \in P_\kappa (\kappa)$, put $\rho_b = \vert Q_b \vert$ and let $\langle Q^{i}_b : i < \rho_b \rangle$ be a one-to-one enumeration of $Q_b$. Now for $\beta < \kappa$ and $i < \rho_\beta$, set $W^{i}_\beta = Q^{i}_\beta \cap \kappa$. We may find a cofinal subset $B$ of $\kappa$ and $h \in \prod_{\beta \in \kappa} \rho_\beta$ such that for any $\alpha, \beta \in B$, we have $W^{h(\alpha)}_\alpha \cap W^{h(\beta)}_\beta \not= \emptyset$. Then clearly, $B \in I_{\kappa, \kappa}^+$, and moreover $Q^{h(\alpha)}_\alpha \cap Q^{h(\beta)}_\beta \not= \emptyset$ for all $\alpha, \beta \in B$.
\hfill$\square$

\medskip

\hskip0,2cm  (ii) $\rightarrow$ (i) :  By Observation 2.7.
\hfill$\square$ 

\begin{Obs} 
 \begin{enumerate} [\rm (i)]
 \item $PS^\ast (\tau, \kappa)$ implies $PS^\ast (\tau, \kappa, \kappa)$.
 \item $PS (\tau, \kappa)$ implies $PS (\tau, \kappa, \kappa)$.
 \end{enumerate}
\end{Obs} 

{\bf Proof.} Argue as for Observation 3.2.
\hfill$\square$

\begin{Cor} $PS (\kappa, \kappa)$ implies $TP (\kappa, \kappa)$.
\end{Cor}

{\bf Proof.} Use Fact 2.22.
\hfill$\square$

\begin{Pro}   \begin{enumerate} [\rm (i)]
\item Suppose that $\lambda$ is regular and $PS^+ (\tau, \kappa, \lambda)$ holds. Then $PS^+ (\tau, \lambda)$ holds.
\item Suppose that $\cf (\lambda) < \kappa$ and $PS^+ ((\cf (\lambda))^+, \kappa, \lambda)$ holds. Then $PS^+ ((\cf (\lambda))^+, \lambda^+)$ holds.
\end{enumerate}
\end{Pro} 

{\bf Proof.} (i) : For $\alpha < \lambda$, let $Q_\alpha$ be a partition of $\lambda \setminus \alpha$ into less than $\tau$ many pieces. For $\alpha \in \lambda$, put $\rho_\alpha = \vert Q_\alpha \vert$ and let $\langle Q^{i}_\alpha : i < \rho_\alpha \rangle$ be a one-to-one enumeration of $Q_\alpha$. Now for $a \in P_\kappa (\lambda)$ and $i < \rho_{\sup a}$, set 

\centerline{$W^{i}_a = \{c \in P_\kappa (\lambda) : a \subseteq c$ and $\sup c \in Q^{i}_{\sup a} \}$.}

 We may find $B \in I_{\kappa, \lambda}^+$ and $h \in \prod_{a \in B} \rho_{\sup a}$ such that for any $a, b \in B$, we have $W^{h(a)}_a \cap W^{h(b)}_b \not= \emptyset$. Put $A = \{\sup a : a \in B\}$, and pick $\psi : A \rightarrow B$ so that $\sup (\psi (\alpha)) = \alpha$ for all $\alpha \in A$. Notice that $A \in I_\lambda^+$. Given $\alpha, \beta \in A$, pick $c \in W^{h(\psi (\alpha))}_{\psi (\alpha)}\cap W^{h(\psi (\beta))}_{\psi (\beta)}$. Then clearly, $\sup c \in Q^{h(\psi (\alpha))}_\alpha \cap Q^{h(\psi (\beta))}_\beta$. 

(ii) : Put $\cf (\lambda) = \sigma$. For $\xi < \lambda^+$, let $W_\xi$ be a partition of $\lambda^+ \setminus \xi$ into at most $\sigma$ many pieces. For $\xi \in \lambda^+$, put $\rho_\xi = \vert W_\xi \vert$ and let $\langle W^{\delta}_\xi : \delta < \rho_\xi \rangle$ be a one-to-one enumeration of $W_\xi$. By Fact 2.10, we may find an increasing sequence $\langle \lambda_i : i < \sigma \rangle$ of regular cardinals greater than $\kappa$ with supremum $\lambda$, and an increasing cofinal sequence $\langle f_\alpha : \alpha < \lambda^+ \rangle$ in $(\prod_{i < \sigma} \lambda_i, <^{\ast})$, where $f <^\ast g$ just in case $\vert \{i < \sigma : f (i) \geq g(i)\} \vert < \sigma$. For $a \in P_\kappa (\lambda)$, define $\chi_a \in \prod_{i < \sigma} \lambda_i$ by $\chi_a (i) = \sup (a \cap \lambda_i)$. Define $s : P_\kappa (\lambda) \rightarrow \lambda^+$ by : $s (a) =$ the least $\alpha$ such that $\chi_a <^\ast f_\alpha$. For $a \in P_\kappa (\lambda)$ and $\delta < \sigma$, let $Q_a^\delta$ denote the collection of all $c \in P_\kappa (\lambda)$ such that $a \subseteq c$ and $s (c) \in W_{s (a)}^\delta$. We may find $B \in I_{\kappa, \lambda}^+$ and $h \in \prod_{a \in B} \rho_{s (a)}$ such that for any $a, b \in B$, we have $Q^{h(a)}_a \cap Q^{h(b)}_b \not= \emptyset$. 

We inductively define $a_n \in B$ for $n < \lambda^+$ so that $s (a_m) < s (a_n)$ whenever $m < n < \lambda^+$. Suppose that $a_m$ has been defined for each $m < n$. Putting $\gamma = \sup \{s (a_m) : m < n \}$, we select $a_n \in B$ so that $ran (f_{\gamma + 1}) \subseteq a_n$. Now let $m < n < \lambda^+$ be given. There must be some $c$ in $Q^{h(a_m)}_{a_m} \cap Q^{h(a_n)}_{a_n}$. Then clearly, $s (c) \geq \max \{ s (a_m), s (a_n)\}$, and moreover $s (c) \in W_{s (a_m)}^{h (a_m)} \cap W_{s (a_n)}^{h (a_n)}$.
\hfill$\square$

\begin{Def} Given an infinite cardinal $\chi$, the \emph{Almost Disjoint Set principle} $ADS_\chi$ asserts the existence of a cofinal subset $y_\alpha$ of $\chi$ of order-type $\cf (\chi)$ for each $\alpha < \chi^+$ such that for each nonzero $\beta < \chi^+$, there is $k \in \prod_{\alpha < \beta} y_\alpha$ with the property that $(y_\delta \setminus (k (\delta)) \cap (y_\alpha \setminus (k (\alpha)) = \emptyset$ whenever $\delta < \alpha < \beta$.
\end{Def}

\medskip

It is known \cite{Scales} that if there is a remarkably good scale on $\chi$, then $ADS_\chi$ holds. Thus the following is closely related to Proposition 2.15.  

\medskip

\begin{Pro}   Let $\chi$ be a singular cardinal such that $PS^+ ((\cf (\chi))^+, \chi)$ holds. Then $ADS_\chi$ fails.
\end{Pro} 

{\bf Proof.} Let $y_\alpha$ be a cofinal subset of $\chi$ of order-type $\cf (\chi)$ for each $\alpha < \chi^+$. Suppose that for each nonzero $\beta < \chi^+$, there is $k_\beta \in \prod_{\alpha < \beta} y_\alpha$ with the property that $(y_\delta \setminus (k_\beta (\delta)) \cap (y_\alpha \setminus (k_\beta (\alpha)) = \emptyset$ whenever $\delta < \alpha < \beta$. For $\alpha < \chi^+$ and $\xi \in y_\alpha$, let $Q^\xi_\alpha$ denote the set of all $\gamma$ such that $\alpha \leq \gamma < \chi^+$ and $k_{\gamma + 1} (\alpha) = \xi$. We may find $B \in [\chi^+]^{\chi^+}$ and $h \in \prod_{\alpha \in B} y_\alpha$ such that $Q^{h (\delta)}_\delta \cap Q^{h (\alpha)}_\alpha \not= \emptyset$ whenever $\delta, \alpha \in B$. Now given $\delta < \alpha < \chi^+$, pick $\gamma \in Q^{h (\delta)}_\delta \cap Q^{h (\alpha)}_\alpha$. Then clearly, $k_{\gamma + 1} (\delta) = h (\delta)$ and $k_{\gamma + 1} (\alpha) = h (\alpha)$, and consequently $(y_\delta \setminus (h (\delta)) \cap (y_\alpha \setminus (h (\alpha)) = \emptyset$. Contradiction.   
\hfill$\square$

\medskip

Let us now turn to the tree property.

\medskip  

%\begin{Obs} $SNS_{\kappa, \kappa} = NS_{\kappa, \kappa}$.\end{Obs}{\bf Proof.} Select a one-to-one function $j$ from $\kappa \times \kappa$ to $\kappa$, and let $D$ be the set of all $a \in P_\kappa (\kappa)$ such that $j (\alpha, \beta) \in a$ whenever $\alpha < \beta$ are in $a$. Now let $C$ be the set of all $a \in P_\kappa (\lambda)$ such that  $j`` ((\gamma + 1) \times (\gamma + 1)) \subseteq a$ for all $\alpha \in a$. Then clearly, $C \subseteq D$. Moreover by Fact 2.20, $C \in SNS_{\kappa,\lambda}^\ast$ and therefore $NS_{\kappa, \lambda} = SNS_{\kappa, \lambda} \vert D = SNS_{\kappa, \lambda}$.\hfill$\square$\begin{Cor} $PS (\kappa, \kappa)$ implies $TP (\kappa, \kappa)$.\end{Cor}{\bf Proof.} By Observation 3.7, $TP (\kappa, \kappa)$ and $TP^- (\kappa, \kappa)$ are equivalent. Now appeal to Fact 2.23 (i) and Observation 3.3 (ii).\hfill$\square$

\begin{Def} The \emph{tree property} $TP (\kappa)$ asserts that any tree of height $\kappa$ each of whose levels has size less than $\kappa$ has a $\kappa$-branch.
\end{Def}

\begin{fact}  \begin{enumerate} [\rm (i)]
\item {\rm(\cite{Weiss})} $TP (\kappa)$ and $TP (\kappa, \kappa)$ are equivalent.
\item {\rm(\cite{Specker})} Let $\tau$ be an infinite cardinal such that $\tau^{< \tau} = \tau$. Then $TP (\tau^+)$ fails.
\end{enumerate}
\end{fact} 

\medskip

 $TP (\kappa)$ can be recast as a piece selection principle.

\medskip
 
 \begin{Obs} The following are equivalent:
 \begin{enumerate} [\rm (i)]
 \item $TP (\kappa)$.
 \item Let $\langle Q_\alpha : \alpha < \kappa \rangle$ be a sequence of partitions of $\kappa$ into less than $\kappa$ many pieces with the property that $Q_\beta \subseteq \bigcup_{W \in Q_\alpha} P (W)$ whenever $\alpha < \beta < \kappa$. Then there is $h \in \prod_{\alpha < \kappa} Q_\alpha$ such that $\vert h (\alpha) \cap h (\beta) \vert \geq 2$ whenever $\alpha < \beta < \kappa$.
 \end{enumerate}
\end{Obs} 

{\bf Proof.}  \hskip0,4cm  (i) $\rightarrow$ (ii) : Suppose that (i) holds, and let $\langle Q_\alpha : \alpha < \kappa \rangle$ be as in (ii). Consider the tree $(T, <_T)$, where
\begin{itemize}
\item $T = \bigcup_{\gamma < \kappa} L_\gamma$, where $L_\gamma$ consists of all $g \in \prod_{\alpha < \gamma + 1} \{A \in Q_\alpha : \vert A \vert \geq 2\}$ such that  $g (\beta) \subseteq g (\alpha)$ whenever $\alpha < \beta < \gamma$.
\item $f <_T g$ just in case $f \subset g$.
\end{itemize}

\hskip0,4cm  (ii) $\rightarrow$ (i) : Suppose that (ii) holds, and let $T = (\kappa, <_T)$ be a tree of height $\kappa$ with each level $L_\alpha$ of size less than $\kappa$. Consider the sequence $\langle Q_\alpha : \alpha < \kappa \rangle$ of partitions of $\kappa$ defined by : $Q_\alpha = \{ Q^\xi_\alpha : \xi \in L_\alpha \}$, where

\centerline{$Q^\xi_\alpha = \{ \zeta \in \kappa : \zeta <_T \xi \} \cup \{ \xi \} \cup \{ \eta \in \kappa : \xi <_T \eta \}$.}
\hfill$\square$

\medskip

This can be used to reformulate $TP (\kappa)$ in terms of partitions relations.

\medskip

 \begin{Obs} The following are equivalent:
 \begin{enumerate} [\rm (i)]
 \item $TP (\kappa)$.
 \item Suppose that $F : \kappa \times \kappa \rightarrow \kappa$ has the following property : if $\beta < \gamma < \delta < \kappa$ are such that $F (\beta, \gamma) = F (\beta, \delta)$, then $F (\alpha, \gamma) = F (\alpha, \delta)$ for all $\alpha < \beta$. Then there is $A \in [\kappa]^\kappa$ such that one of the following holds :
 \begin{itemize}
 \item $F (\beta, \gamma) \not= F (\beta, \delta)$ whenever $\beta < \gamma < \delta$ are in $A$.
 \item $F (\beta, \gamma) = F (\beta, \delta)$ whenever $\beta < \gamma < \delta$ are in $A$.
 \end{itemize}
  \end{enumerate}
\end{Obs} 

{\bf Proof.}  \hskip0,4cm  (i) $\rightarrow$ (ii) : Assume that (i) holds, and let $F$ be as in (ii). For $\alpha < \kappa$, consider the equivalence relation $\thicksim_\alpha$ defined on $\kappa$ by : $\beta \thicksim_\alpha \gamma$ if and only if either $\gamma = \beta \leq \alpha$, or $\beta, \gamma > \alpha$ and $F (\xi, \beta) = F (\xi, \gamma)$ for all $\xi \leq \alpha$. Let $Q_\alpha$ be the set of all equivalence classes with respect to $\thicksim_\alpha$.

\underline {Case 1} : There is $\eta < \kappa$ such that $\vert Q_\eta \vert = \kappa$. Pick $A \in [\kappa \setminus (\eta + 1))]^\kappa$ so that $\vert A \cap H \vert \leq 1$ for all $H \in Q_\eta$. Now if $\beta < \gamma < \delta$ are in $A$, we must have $F (\beta, \gamma) \not= F (\beta, \delta)$, since otherwise we would have $F (\eta, \gamma) = F (\eta, \delta)$.

\underline {Case 2} : $\vert Q_\alpha \vert < \kappa$ for all $\alpha < \kappa$. Then by Observation 3.10, we may find $h \in \prod_{\alpha < \kappa} Q_\alpha$ such that $\vert h (\alpha) \cap h (\beta) \vert \geq 2$ whenever $\alpha < \beta < \kappa$. It is simple to see that $\{h (\alpha) : \alpha < \kappa\} \subseteq [\kappa]^\kappa$. Furthermore, $h (\beta) \subseteq h (\alpha)$ whenever $\alpha < \beta < \kappa$. Now inductively define an increasing sequence $\langle \alpha_i : i < \kappa \rangle$ of elements of $\kappa$ so that for any $j < \kappa$, $\alpha_j \in h ((\sup \{\alpha_i : i < j\}) + 1)$. Then clearly, $F (\alpha_i, \alpha_j) = F (\alpha_i, \alpha_k)$ whenever $i < j < k < \kappa$.

\hskip0,4cm  (ii) $\rightarrow$ (i) : Assume that (ii) holds, and let $\langle Q_\alpha : \alpha < \kappa\rangle$ be as in (ii) of Observation 3.10. For $\alpha < \kappa$, let $\langle Q^{i}_\alpha : i < \vert Q_\alphaÊ\vert \rangle$ be a one-to-one enumeration of $Q_\alpha$. Define $F : \kappa \times \kappa \rightarrow \kappa$ by : $F (\alpha, \beta) = i$ just in case $\beta \in Q^{i}_\alpha$. There must be $A \in [\kappa]^\kappa$ and $f \in \prod_{\alpha \in A} \vert Q_\alpha \vert$ such that $F (\beta, \gamma) = f (\beta)$ whenever $\beta < \gamma$ are in $A$. Then $A \setminus (\beta + 1) \subseteq Q^{f (\beta)}_\beta$ for all $\beta \in A$. It easily follows that the conclusion of (ii) of Observation 3.10 holds.
\hfill$\square$

\medskip

{\bf QUESTION.} Is it consistent that  $TP (\kappa)$ holds, but $PS (\kappa, \kappa)$ fails ?

\medskip

We return to ideal extension, but this time for ideals on $\kappa$.

\medskip

\begin{Def} We let $L_{\kappa\omega}$ denote the infinitary language which allows conjunctions and disjunctions of less than $\kappa$ many formulas, and universal and existential quantification over finitely many variables.

$L_{\kappa\omega}$ is \emph{weakly compact} if any set of $\kappa$ sentences from $L_{\kappa\omega}$ without a model has a subset of smaller size without a model.
\end{Def}

\begin{Def}  For a fine ideal $J$ on $\kappa$, $IE (\tau, \kappa, J)$ means that given a partition $Q_\alpha$ of $\kappa \setminus \alpha$ with $\vert Q_\alpha \vert < \tau$ for each $\alpha \in \kappa$, there is $h \in \prod_{\alpha \in \kappa} Q_\alpha$ and an ideal $K$ on $\kappa$ extending $J$ such that $ran (h) \subseteq K^\ast$.
\end{Def}

\begin{Obs} Suppose that $L_{\kappa\omega}$ is weakly compact. Then $IE (\kappa, \kappa, I_\kappa)$ holds.
\end{Obs}

{\bf Proof.} For $\alpha < \kappa$, let $Q_\alpha$ be a partition of $\kappa \setminus \alpha$ with $\vert Q_\alpha \vert < \kappa$. Consider the $L_{\kappa\omega}$ language with one unary predicate $S$ and constant symbols $c_A$ for $A \in \bigcup_{\alpha < \kappa} Q_\alpha$. Let $\Sigma$ consist of the following sentences :
\begin{itemize}
\item $\bigvee_{A \in Q_\alpha} S (c_A)$ for each $\alpha < \kappa$.
\item $\neg (S (c_{A_0}) \wedge S (c_{A_1}) \wedge \cdots \wedge S (c_{A_n}))$ whenever $0 < n < \omega$, $A_0, A_1, \cdots, A_n \in  \bigcup_{\alpha < \kappa} Q_\alpha$ and $A_0 \cap A_1 \cap \cdots \cap A_n = \emptyset$.  
\end{itemize} 
Notice that for $0 < \beta < \kappa$, $\bigcap_{\alpha < \beta} k_\beta (\alpha) \not= \emptyset$, where $k_\beta : \beta \rightarrow \kappa$ is defined by $k_\beta (\alpha) =$ the unique $A \in Q_\alpha$ such that $\beta \in A$. It easily follows that any subset of $\Sigma$ of size less than $\kappa$ is satisfiable. Hence so is $\Sigma$ itself, and there must be $h \in \prod_{\alpha < \kappa} Q_\alpha$ with the property that for each $e \in P_\omega (\kappa) \setminus \{\emptyset\}$, $\bigcap_{\alpha \in e} h (\alpha)$ is nonempty. In fact,  $\bigcap_{\alpha \in e} h (\alpha) \in I_\kappa^+$. Suppose otherwise, and let $\delta < \kappa$ such that $\bigcap_{\alpha \in e} h (\alpha) \subseteq \delta$. Then $\bigcap_{\alpha \in d} h (\alpha) = \emptyset$, where $d = e \cup \{ \delta \}$. Contradiction.
\hfill$\square$

\medskip

Boos \cite{Boos} showed that if $\kappa$ is weakly compact, then in the extension obtained by adding $\kappa^+$ many Cohen reals, $L_{\kappa\omega}$ is still weakly compact. Thus $IE (\kappa, \kappa, I_\kappa)$ (and hence $PS^+ (\kappa, \kappa)$) may hold without $\kappa$ being inaccessible.

\medskip

{\bf QUESTION.} Is it consistent that  $PS^+ (\kappa, \lambda')$ holds for every cardinal $\lambda' \geq \kappa$, but $\kappa$ is not inaccessible ?

\medskip 

\begin{Def} For an infinite cardinal $\tau$, the \emph{transversal property} $PT (\kappa, \tau)$ means that for any size $\kappa$ family of sets of size less than $\tau$ without a transversal (i.e. a one-to-one choice function), there exists a subfamily of size less than $\kappa$ without a transversal.
\end{Def}

\begin{fact}  {\rm(\cite{MS})}
%\begin{enumerate}[(i)]\item Let $\pi < \omega_{\omega^2+1}$ be a regular infinite cardinal. Then {\rm PT}$(\pi, \omega_1)$ does not hold.\item It is consistent (relative to infinitely many supercompact cardinals) that {\rm PT}$(\omega_{\omega^2+1}, \omega_1)$ holds\item Let $\pi$ be a regular infinite cardinal such that {\rm PT}$(\pi,\omega_1)$ does not hold. Then {\rm PT}$(\omega_{\pi+1}, \omega_1)$ does not hold.\item 
It is consistent (relative to infinitely many supercompact cardinals) that $PT (k, \omega_1)$ holds for every regular infinite cardinal $k$ greater than the least fixed point of the aleph function.
%\end{enumerate}
\end{fact} 

\begin{Obs} Suppose that $IE (\tau, \kappa, I_\kappa)$ holds. Then so does $PT (\kappa, \tau)$.
\end{Obs}

{\bf Proof.} Let $\langle X_\alpha : \alpha < \kappa \rangle$ be a sequence of sets of size less than $\tau$ with the property that for any nonzero $\beta < \kappa$, there is a one-to-one $k_\beta$ in $\prod_{\alpha < \beta} X_\alpha$. Pick a one-to-one function $j : \bigcup_{\alpha < \kappa} X_\alpha \rightarrow \kappa$. For $\alpha < \kappa$ and $i \in  j``X_\alpha$, set $A^{i}_\alpha = \{ \gamma \in \kappa \setminus \alpha : j (k_{\gamma + 1} (\alpha)) = i \}$. There must be $h \in \prod_{\alpha < \kappa} j``X_\alpha$ such that $A^{h (\alpha)}_\alpha \cap A^{h (\beta)}_\beta \not= \emptyset$ for all $\alpha, \beta < \kappa$. Define $g \in \prod_{\alpha < \kappa} X_\alpha$ so that $j (g (\alpha)) = h (\alpha)$. We will show that $g$ is one-to-one. Thus let $\alpha < \beta < \kappa$. Select $\gamma$ in $A^{h (\alpha)}_\alpha \cap A^{h (\beta)}_\beta$. Then $ j (k_{\gamma + 1} (\alpha)) = h (\alpha)$ and $ j (k_{\gamma + 1} (\beta)) = h (\beta)$, which gives $k_{\gamma + 1} (\alpha) = g (\alpha)$ and $k_{\gamma + 1} (\beta) = g (\beta)$. It follows that $g (\alpha) \not= g (\beta)$.
\hfill$\square$

\medskip 

Let $PT^- (\kappa, \tau)$ mean that for any size $\kappa$ family of sets of size less than $\tau$ with the property that any subfamily of size less than $\kappa$ has a  transversal, there exists a subfamily of size $\kappa$ with a transversal. Then by the proof of Observation 3.17, $PS^+ (\tau, \kappa, \tau)$ implies $PT^- (\kappa, \tau)$.
 
\medskip

{\bf QUESTION.} Is it consistent that  $PS^+ (\kappa, \kappa)$ holds, but $IE (\kappa, \kappa, I_\kappa)$  fails ?

\medskip

{\bf QUESTION.} What is the least possible value of $\kappa$ at which $PS^+ (\kappa, \kappa)$ (respectively, $PS^\ast (\kappa, \kappa)$, $PS (\kappa, \kappa)$) may hold ? 
 
\bigskip

\section{Covering numbers}
  
\medskip
 
In this section we study the consequences of $PS^+$ in terms of cardinal arithmetic (in the sense of Shelah).

\medskip  

\begin{Def} Given two infinite cardinals $\rho \leq \sigma$, $u(\rho, \sigma)$ denotes the cofinality of the poset $(P_\rho (\sigma), \subseteq)$.  
 \end{Def}

\begin{fact} {\rm (Folklore)} Let $\rho \leq \sigma$ be two infinite cardinals. Then $\sigma^{< \rho} = \max\{2^{<\rho}, u(\rho, \sigma)\}$.
\end{fact}

%\begin{fact}    {\rm(\cite{DM})}  \begin{enumerate}[\rm (i)]\item  Let $\pi$ be a regular cardinal such that $\kappa \leq \pi \leq \lambda$. Then $u(\kappa, \lambda) \leq \max\{u(\kappa, \pi), u(\pi, \lambda)\}$.  \item  Let $\sigma$ be a cardinal such that $\kappa \leq \sigma < \lambda$. Then $u(\kappa, \lambda) \leq \max\{u(\kappa, \sigma), u(\sigma^+, \lambda)\}$.   \end{enumerate}\end{fact}

 \begin{Def}  Given four cardinals $\rho_1, \rho_2, \rho_3, \rho_4$ with $\rho_1 \geq \rho_2 \geq \rho_3 \geq \omega$ and $\rho_3 \geq \rho_4 \geq 2$, $\cov (\rho_1, \rho_2, \rho_3, \rho_4)$  denotes the least cardinality of any $X \subseteq  P_{\rho_2}(\rho_1)$ such that for any $a \in  P_{\rho_3}(\rho_1)$, there is $Q \in  P_{\rho_4}(X)$ with $a \subseteq \bigcup Q$.          
\end{Def}

\medskip

Note that $u(\rho, \sigma) = \cov (\sigma, \rho, \rho, 2)$.  

\medskip 

\begin{fact}   {\rm(\cite[pp. 85-86]{SheCA}, \cite{LCCN})}  Let $\rho_1, \rho_2, \rho_3$ and $\rho_4$ be four cardinals such that $\rho_1 \geq \rho_2 \geq \rho_3 \geq \omega$ and $\rho_3 \geq \rho_4 \geq 2$. Then the following hold :        
\begin{enumerate}[\rm (i)]
\item  If $\rho_1 = \rho_2$ and either $\cf(\rho_1) < \rho_4$ or $\cf(\rho_1) \geq \rho_3$, then $\cov (\rho_1, \rho_2, \rho_3, \rho_4) = \cf(\rho_1)$.           
\item  If either $\rho_1 > \rho_2$, or $\rho_1 = \rho_2$ and $\rho_4 \leq \cf(\rho_1) < \rho_3$, then $\cov (\rho_1, \rho_2, \rho_3, \rho_4) \geq \rho_1$.        
\item  $\cov (\rho_1, \rho_2, \rho_3, \rho_4) = \cov (\rho_1, \rho_2, \rho_3, \max\{\omega, \rho_4\})$.
\item $\cov (\rho_1^+, \rho_2, \rho_3, \rho_4) = \max\{\rho_1^+, \cov (\rho_1, \rho_2, \rho_3, \rho_4)\}$.
\item  If $\rho_1 > \rho_2$ and $\cf(\rho_1) < \rho_4 = \cf(\rho_4)$, then 

\centerline{$\cov (\rho_1, \rho_2, \rho_3, \rho_4) = \sup \{\cov (\rho, \rho_2, \rho_3, \rho_4) : \rho_2 \leq \rho < \rho_1\}$.}

\item   If $\rho_1$ is a limit cardinal such that $\rho_1 > \rho_2$ and $\cf(\rho_1) \geq \rho_3$, then 

\centerline{$\cov (\rho_1, \rho_2, \rho_3, \rho_4) = \sup\{\cov (\rho, \rho_2, \rho_3, \rho_4) : \rho_2 \leq \rho < \rho_1\}$.}

\item   If $\rho_3 > \rho_4 \geq \omega$, then
 
\centerline{$\cov (\rho_1, \rho_2, \rho_3, \rho_4) = \sup\{\cov (\rho_1, \rho_2, \rho^+, \rho_4) : \rho_4 \leq \rho < \rho_3\}$.}
       
\item   If $\rho_3 \leq \rho_2 = \cf(\rho_2)$, $\omega \leq \rho_4 = \cf(\rho_4)$ and $\rho_1 < \rho_2^{+\rho_4}$, then $\cov (\rho_1, \rho_2, \rho_3, \rho_4) = \rho_1$. 
\item   If $\rho_3 = \cf(\rho_3)$, then either  $\cf(\cov (\rho_1, \rho_2, \rho_3, \rho_4)) < \rho_4$, or  $\cf(\cov (\rho_1, \rho_2, \rho_3, \rho_4)) \geq \rho_3$.
\item Suppose that $\rho_3 > \cf (\rho_2) \geq \rho_4$ and $\cf (\rho_3) \not= \cf (\rho_2)$. Then $\cov (\rho_1, \rho_2, \rho_3, \rho_4) = \cov (\rho_1, \rho, \rho_3, \rho_4)$ for some cardinal $\rho$ with $\rho_2 > \rho \geq \rho_3$.
 \end{enumerate}
\end{fact}

\begin{fact}  \begin{enumerate}[\rm (i)]
\item {\rm(\cite[Remark 6.6.A  p. 101]{SheCA})} Let $\chi$ be a singular cardinal. Then $\cov (\chi,\chi, (\cf (\chi))^+, \cf (\chi)) > \chi^+$ if and only if $\cov (\chi,\chi, (\cf (\chi))^+, 2) > \chi^+$.
\item  {\rm(\cite[p. 99]{SheCA}, \cite{LCCN}))} Let $\chi$ be a singular cardinal. Suppose that $\cov (\chi, \chi, (cf (\chi))^+, 2) > \chi^+$. Then we may find $y_\alpha \in P_{(cf (\chi))^+} (\chi)$ for $\alpha < \chi^+$ such that for any nonzero $\beta < \chi^+$, there is a one-to-one $h \in \prod_{\alpha < \beta} y_\alpha$.
\end{enumerate}
\end{fact}
 
%{\bf Proof.} Select a partition $\langle X_\delta : \delta < \tau \rangle$ of $\chi$ into $\tau$ pieces of size $\chi$. For $\beta \in \chi^+ \setminus \{0\}$, let $\beta = \bigcup_{\gamma < \cf (\chi)} d^\beta_\gamma$, where $\vert d_\gamma^\beta \vert < \chi$ for all $\gamma < \cf (\chi)$. Define $g : \bigcup_{0 < \beta < \chi^+} (\{\beta\} \times \beta) \rightarrow \cf (\chi)$ by $g (\beta, \alpha) =$ the least $\gamma < \cf (\chi)$ such that $\alpha \in d^\beta_\gamma$.Given $\delta < \chi$, define by induction $y_\alpha^\delta \in P_{\tau^+}(X_\delta)$ for $\alpha < \chi^+$ so that for every $\beta < \chi^+$,\centerline{$y_\alpha^\delta \setminus (\bigcup \{ y_\zeta^\delta : \zeta < \alpha$ and  $g (\beta, \zeta) = g (\beta, \alpha)\}) \not= \emptyset$.}For $\alpha < \chi^+$, set $y_\alpha = \bigcup_{\delta < \tau} y_\alpha^\delta$. Now suppose $0 < \beta < \chi^+$.  Select\centerline{$h \in \prod_{\alpha < \beta} y_\alpha^\delta \setminus (\bigcup \{ y_\zeta^\delta : \zeta < \alpha$ and  $g (\beta, \zeta) = g (\beta, \alpha)\})$.}Then $h$ is easily seen to be one-to-one.  \hfill$\square$
 
\medskip

Note that if $\vec{f} = \langle f_\alpha : \alpha < \pi \rangle$ is an increasing, cofinal sequence in $(\prod A, <_I)$, then by a result of Shelah \cite[Theorem 5.4 pp. 87-88]{SheCA}, $\pi \leq \cov (\sup A, \sup A, \vert A \vert^+, 2)$.

\medskip

 \begin{Obs}  Let $\chi$ be a singular cardinal such that $PS^+ ((\cf (\chi))^+, \chi^+)$ holds. Then  $\cov (\chi, \chi, (\cf (\chi))^+, 2) = \chi^+$.
 \end{Obs}
 
{\bf Proof.}   Suppose otherwise. Put $\sigma = \cf (\chi)$. By Fact 4.5, we may find $y_\alpha \in P_{\sigma^+} (\chi)$ for $\alpha < \chi^+$ such that for any nonzero $\beta < \chi^+$, there is a one-to-one $h_\beta \in \prod_{\alpha < \beta} y_\alpha$. For $\alpha < \chi^+$, let $\langle y^{i}_\alpha : i < \vert y_\alpha \vert \rangle$ be a one-to-one enumeration of $y_\alpha$, and set  $Q^{i}_\alpha = \{\delta \in \chi^+ \setminus \alpha : h_{\delta + 1} (\alpha) = y^{i}_\alpha\}$. There must be $B \in [\chi^+]^{\chi^+}$ and $g \in \prod_{\alpha \in B} \vert y_\alpha \vert$ with the property that $Q^{g (\alpha)}_\alpha \cap Q^{g (\gamma)}_\gamma \not= \emptyset$ whenever $\alpha, \gamma \in B$. Then clearly, the function $t : \chi^+ \rightarrow \chi$ defined by $t (\alpha) = y_\alpha^{g (\alpha)}$ is one-to-one. Contradiction.
%Pick a bijection $j : \lambda^+ \times \omega \rightarrow \lambda^+$. For $a \in P_\kappa (\lambda^+)$, let \centerline{$s_a = a \cap \{j (\alpha, n) : \alpha < \sup a$ and $h_{\sup a} (\alpha) = y^n_\alpha\}$.}Define $\psi : P_3 (\lambda^+) \setminus \{\emptyset\} \rightarrow P_{\omega_1} (\lambda^+)$ by $\psi (e) = e \cup \{(\max e) + 1\} \cup \{j (\alpha, n) : \alpha \in e$. There must be $F \subseteq \lambda^+ \times \omega$ with the property that for any $ e \in P_3 (\lambda) \setminus \{\emptyset\}$, there is $a \in P_\kappa (\lambda^+)$ such that $\psi (e) \subseteq a$ and $j`` F \cap \psi (e) = s_a \cap \psi (e)$. Then it is easy to see that $F$ is a function with domain $\lambda^+$. Moreover, the function $g : \lambda^+ \rightarrow \lambda$ defined by $g (\alpha) = y_\alpha^{F (\alpha)}$ is one-to-one. Contradiction.
 \hfill$\square$

\medskip

Neeman \cite{Itai} established the consistency relative to large cardinals of the existence of a singular strong limit cardinal $\chi$ of cofinality $\omega$ such that $TP (\chi^+)$ holds and $2^\chi > \chi^+$. Note that

\centerline{$2^\chi = \chi^{\cf (\chi)} \leq \max \{\cov (\chi, \chi, (\cf (\chi))^+, 2), 2^{< \chi}\} = \cov (\chi, \chi, (\cf (\chi))^+, 2) \leq 2^\chi$,}

and therefore by Observation 4.6, $PS^+ ((\cf (\chi))^+, \chi^+)$ fails. In Neeman's model, there is both a very good (and hence remarkably good) scale of length $\chi^+$ on $\chi$, and a scale of length $\chi^+$ on $\chi$ that is not good (so that approachability fails, and in fact \cite{Scales} there is a regular uncountable cardinal $\sigma < \chi$ such that $E^{\chi^+}_\sigma \notin I[\chi^+ ; \chi]$).

\medskip

To get the most out of Observation 4.6 we will vary the value of $\lambda$. We will thus be able to use the following results of pcf theory.

\medskip
  
 \begin{fact} {\rm (\cite{LCCN})} Let $\sigma$, $k$, $\mu$ and $\nu$ be four infinite cardinals such that $\cf (\sigma) = \sigma \leq k$ and $\sigma < \cf (\mu) = \mu < \nu$. Suppose that 
 \begin{itemize}
 \item $\cov (k, \rho^+, \rho^+, \sigma) \leq k^{++}$ for every cardinal $\rho$ with $\sigma \leq \rho < \min \{k, \mu\}$.
 \item  $\cov (\chi, \chi, \sigma^+, \sigma) = \chi^+$ for every cardinal $\chi$ with $k < \chi \leq \nu$ and $\cf (\chi) = \sigma$.
\end{itemize}
Then $\cov (\nu, \mu, \mu, \sigma) \leq \nu^+$.
 \end{fact}

 \begin{Obs}
 \begin{enumerate} [\rm (i)]
 \item Let $\pi$ and $\mu$ be two regular cardinals such that $\omega \leq \pi \leq \kappa \leq \mu \leq \lambda$. Suppose that 
 \begin{itemize}
 \item either $\cf (\lambda) < \pi$, or $\cf (\lambda) \geq \mu$.
  \item $u (\rho^+, \kappa) \leq \kappa^{++}$ for every cardinal $\rho$ with $\omega \leq \rho < \kappa$.
 \item  $\cov (\chi, \chi, \omega_1, 2) = \chi^+$ for every cardinal $\chi$ with $\kappa < \chi < \lambda$ and $\cf (\chi) = \omega$.
\end{itemize}
Then $\cov (\lambda, \mu, \mu, \pi) = \lambda$.
\item Let $\pi$ and $\mu$ be two regular cardinals such that $\omega_1 \leq \pi \leq \kappa \leq \mu \leq \lambda$. Suppose that 
 \begin{itemize}
 \item either $\cf (\lambda) < \pi$, or $\cf (\lambda) \geq \mu$.
  \item $\cov (\kappa, \rho^+, \rho^+, \omega_1) \leq \kappa^{++}$ for every cardinal $\rho$ with $\omega_1 \leq \rho < \kappa$.
 \item  $\cov (\chi, \chi, \omega_2, \omega_1) = \chi^+$ for every cardinal $\chi$ with $\kappa < \chi < \lambda$ and $\cf (\chi) = \omega_1$.
\end{itemize}
Then $\cov (\lambda, \mu, \mu, \pi) = \lambda$.
 \end{enumerate}
\end{Obs} 

{\bf Proof.}  We prove (i) and leave the similar proof of (ii) to the reader. By Fact 4.4 ((i) and (ii)), $\cov (\tau, \mu, \mu, \pi) \geq \tau$ for every cardinal $\tau \geq \mu$. Furthermore by Fact 4.7, $\cov (\nu, \mu, \mu, \pi) \leq u (\mu, \nu) \leq \nu^+$ for any cardinal $\nu$ with $\mu < \nu < \lambda$.

\underline {Case 1} : $\lambda = \mu$. Then by Fact 4.4 (i),

\centerline{$\lambda \leq \cov (\lambda, \mu, \mu, \pi) \leq u (\mu, \lambda) = \lambda$.}

\underline {Case 2} : $\lambda$ is the successor of some cardinal $\sigma \geq \mu$. Then by Fact 4.4 (i), 

\centerline{$\lambda \leq \cov (\lambda, \mu, \mu, \pi) \leq u (\mu, \lambda) = \max \{\lambda, u (\mu, \sigma)\} = \lambda$.}

\underline {Case 3} : $\cf (\lambda) < \pi$. Then by Fact 4.4 (v), 

\centerline{$\lambda \leq \cov (\lambda, \mu, \mu, \pi) = \sup \{\cov (\nu, \mu, \mu, \pi) : \mu \leq \nu < \lambda\} \leq \lambda$.}

\underline {Case 4} : $\lambda$ is a limit cardinal with $\mu \leq \cf (\lambda)$. Then by Fact 4.4 (vi), 

\centerline{$\lambda \leq \cov (\lambda, \mu, \mu, \pi) = \sup\{\cov (\nu, \mu, \mu, \pi) : \mu \leq \nu < \lambda\} \leq \lambda$.}
\hfill$\square$

\begin{fact}    {\rm(\cite{DM})}  
\begin{enumerate}[\rm (i)]
\item  Let $\pi$ be a regular cardinal such that $\kappa \leq \pi \leq \lambda$. Then $u (\kappa, \lambda) \leq \max \{u (\kappa, \pi), u (\pi, \lambda)\}$.  
\item Suppose that $\lambda$ is a limit cardinal. Then

\centerline{$u (\kappa, \lambda) = \max \{\cov (\lambda, \lambda, \kappa, 2),  \sup \{u (\kappa, \chi) : \kappa \leq \chi < \lambda\}\}$.}

\end{enumerate}
\end{fact}

 \begin{Obs} Suppose that 
 \begin{itemize}
 \item $\cf (\lambda) = \omega$.
  \item $u (\rho^+, \kappa) \leq \kappa^{++}$ for every cardinal $\rho$ with $\omega \leq \rho < \kappa$.
 \item  $\cov (\chi, \chi, \omega_1, 2) = \chi^+$ for every cardinal $\chi$ with $\kappa < \chi < \lambda$ and $\cf (\chi) = \omega$.
\end{itemize}
Then $u (\kappa, \lambda) = \cov (\lambda, \lambda, \omega_1, 2)$.
\end{Obs} 

{\bf Proof.}  \medskip

{\bf Claim 1.}   Let $\chi$ be a cardinal with $\kappa^{++} \leq \chi \leq \lambda$. Then $u (\omega_1, \chi)  \leq u (\kappa, \chi)$. 

\smallskip

{\bf Proof of Claim 1.} By Facts 4.4 and  4.9 (i),  

\centerline{$u (\omega_1, \chi)  \leq \max \{u (\omega_1, \kappa), u (\kappa, \chi)\} = u (\kappa, \chi)$,}

 which completes the proof of the claim.

\medskip

{\bf Claim 2.}   $u (\kappa, \lambda)  = u (\omega_1, \lambda)$. 

\smallskip

{\bf Proof of Claim 2.} By Observation 4.8 (i),  $\cov (\lambda, \lambda, \omega_1, \omega_1) = \lambda$. Hence by Claim 1, 

\centerline{$u (\kappa, \lambda) \leq u (\omega_1, \cov (\lambda, \lambda, \omega_1, \omega_1)) = u (\omega_1, \lambda)  \leq  u (\kappa, \lambda)$,}

 which completes the proof of the claim.

\medskip

{\bf Claim 3.}   $u (\omega_1, \lambda) = \cov (\lambda, \lambda, \omega_1, 2)$.

\smallskip

{\bf Proof of Claim 3.} By Claim 1 and Facts 4.4 and 4.9 (ii), 

\centerline{$u (\omega_1, \lambda) = \max \{\cov (\lambda, \lambda, \omega_1, 2),  \sup \{u (\omega_1, \tau) : \omega_1 \leq \tau < \lambda\}\} = \cov (\lambda, \lambda, \omega_1, 2)$,}

which completes the proof of the claim and that of the observation.
\hfill$\square$

% Let $\sigma$, $\mu$ and $\nu$ be three infinite cardinals such that $\cf (\sigma) = \sigma < \kappa \leq \mu < \nu$. Suppose that \begin{itemize}\item $\cov (\kappa, \rho^+, \rho^+, \sigma) \leq \kappa^{++}$ for every cardinal $\rho$ with $\sigma \leq \rho < \kappa$.\item  $PS (\sigma^+, \kappa, \chi)$ holds for every cardinal $\chi$ with $\kappa < \chi \leq \nu$ and $\cf (\chi) = \sigma$\end{itemize}Then $\cov (\nu, \mu, \mu, \sigma) \leq \nu^+$. \end{Pro}

\bigskip

\section{Unbalanced partition properties}

\medskip

Let us start with partitions of $P_\kappa (\lambda) \times P_\kappa (\lambda)$.

\medskip

\begin{Def} Given two collections $X$ and $Y$ of subsets of $P_\kappa (\lambda)$, an ideal $J$ on $P_\kappa (\lambda)$, and a cardinal $\rho$ with $0 < \rho \leq \kappa$, $X \xrightarrow{Y} (J^+, \rho)^2$ means that for any $F : P_\kappa (\lambda) \times P_\kappa (\lambda) \rightarrow 2$ and any $A \in X$, there is either $B \in J^+ \cap P (A)$ such that $\{ b \in B : F (a, b) \not= 0\} \in Y$ for all $a \in B$, or an increasing sequence $\langle c_i : i < \rho \rangle$ in $(A, \subset)$ such that $F (c_i, c_j) = 1$ whenever $i < j < \rho$.  
\end{Def}

\begin{Obs} Let $J$ be a fine ideal on $P_\kappa (\lambda)$. Then  $J^+ \xrightarrow{J} (J^+, \omega)^2$ holds. 
\end{Obs}

{\bf Proof.}  Fix $F : P_\kappa (\lambda) \times P_\kappa (\lambda) \rightarrow 2$ and $A \in J^+$. Define $\psi : P_\kappa (\lambda) \times P (A) \rightarrow P (A)$ by 

\centerline{$\psi (a, X) = \{b \in X : a \subset b$ and $F (a, b) = 1\}$.}
 
\underline {Case 1} : There is $Y \in J^+  \cap P (A)$ such that $\psi (a, Y) \in J$ for all $a \in Y$. Then clearly, $F (a, b) = 0$ whenever $a$ is in $Y$ and $b$ is in $Y \setminus \psi (a, Y)$ with $a \subset b$.

\underline {Case 2} : For each $X \in J^+ \cap P (A)$, there is $a_X$ in $X$ such that $\psi (a_X, X) \in J^+$. Inductively define $X_n$ for $n < \omega$ by :
\begin{itemize}
\item $X_0 = A$.
\item $X_{n + 1} = \psi (a_{X_n}, X_n)$.
\end{itemize}
Then clearly, $\{a_{X_n} : n < \omega\} \subseteq A$. Moreover, if $m <  n < \omega$, then $a_{X_m} \subset a_{X_n}$ and $F (a_{X_m}, a_{X_n}) = 1$.
\hfill$\square$

\begin{Def} For an ideal $J$ on a set $X$, $MAD (J)$ (respectively, $MAD_d (J)$) denotes the collection of all $Q \subseteq J^+$ such that
\begin{itemize}
\item $A \cap B \in J$ (respectively, $A \cap B = \emptyset$) for any two distinct members $A, B$ of $Q$.
\item For any $C \in J^+$, there is $A \in Q$ with $A \cap C \in J^+$.  
\end{itemize}

Let $\rho$ and $\nu$ be two nonzero cardinals. $J$ is $(\rho, \nu)$-\emph{distributive} (respectively, \emph{disjointly} $(\rho, \nu)$-\emph{distributive}) if given $A \in J^+$, and $Q_\alpha$ in $MAD (J)$ (respectively, $MAD_d (J)$) with $\vert Q_\alpha \vert \leq \nu$ for $\alpha < \rho$, there is $B \in J^+ \cap P (A)$ and $h \in \prod_{\alpha < \rho} Q_\alpha$ such that $B \setminus h (\alpha) \in J$ for every $\alpha < \rho$. 
\end{Def}

\begin{Obs} Suppose that $J$ is a $(\theta, 2)$-distributive ideal on a set $X$, where $\theta$ is an infinite cardinal. Then $J$ is $(\theta, \theta)$-distributive. 
\end{Obs}

{\bf Proof.}  Let $A \in J^+$, and $Q_\alpha \in MAD (J)$ with $\vert Q_\alpha \vert \leq \theta$ for $\alpha < \theta$. Put $Z = \bigcup_{\alpha < \theta} Q_\alpha$. There must be $B \in J^+ \cap P (A)$ and $h \in \prod_{W \in Z} \{W, X \setminus W\}$ such that $B \setminus h (W) \in J$ for every $W \in Z$. Now given  $\alpha < \theta$,  we have $Q_\alpha \in MAD (J)$, so there is a (unique) $W \in Q_\alpha$ such that $h (W) = W$.   
\hfill$\square$
  
\begin{Obs} Suppose that $J^+ \xrightarrow{J} (J^+, \sigma^+)^2$ holds, where $J$ is a $\kappa$-complete, fine ideal on $P_\kappa (\lambda)$, and $\sigma$ is an infinite cardinal. Then $J$ is disjointly $(\sigma, \lambda^{< \kappa})$-distributive. 
\end{Obs}

{\bf Proof.}  Let $A \in J^+$, and $Q_\alpha$ in $MAD_d (J)$ for $\alpha < \sigma$. For $\alpha < \sigma$, put $X^0_\alpha = P_\kappa (\lambda) \setminus \bigcup Q_\alpha$, and let $\langle X^{i}_\alpha : 0 < i < \vert Q_\alpha \vert \rangle$ be a one-to-one enumeration of $Q_\alpha$. Define $g : \sigma \times P_\kappa (\lambda) \rightarrow\lambda^{< \kappa}$ so that for any $a \in P_\kappa (\lambda)$ and any $\alpha < \sigma$, $a \in X^{g (\alpha, a)}_\alpha$. Now define $F : P_\kappa (\lambda) \times P_\kappa (\lambda) \rightarrow 2$ by : $F (a, b) = 1$ if and only if there is $\alpha < \sigma$ such that $g (\alpha, a) \not= g (\alpha, b)$, and for the least such $\alpha$, $g (\alpha, a) > g (\alpha, b)$.

\underline {Case 1} : There is $B \in J^+ $, and $Z_a \in J$ for $a \in B$ such that $F (a, b) = 0$ whenever $a, b \in B$ and $b \notin Z_a$. Assume toward a contradiction that there is $\gamma < \sigma$ such that $\{ B \setminus X^{i}_\gamma : 0 < i < \vert Q_\gamma \vert \} \subseteq J^+$, and let $\alpha$ denote the least such $\gamma$. Let $h \in \prod_{\beta < \alpha} \{k : 0 < k < \vert Q_\beta \vert \}$ such that $\{ B \setminus X^{h (\beta)}_\beta : \beta < \alpha \} \subseteq J$. There must be $0 < i < j < \vert Q_\alpha\vert $ such that $\{B \cap X^{i}_\alpha, B \cap X^j_\alpha \} \subseteq J^+$. Pick $a$ in $(B \cap X^j_\alpha) \setminus \bigcup_{\beta < \alpha} (B \setminus X^{h (\beta)}_\beta)$, and $b$ in $(B \cap X^{i}_\alpha) \setminus (Z_a \cup \bigcup_{\beta < \alpha} (B \setminus X^{h (\beta)}_\beta))$. Then $F (a, b) = 1$. Contradiction.

 \underline {Case 2} : There is an increasing sequence $\langle a_\delta : \delta < \sigma^+ \rangle$ in $(A, \subset)$ such that $F (a_\gamma, a_\delta) = 1$ whenever $\gamma < \delta < \sigma^+$. We will show that this is contradictory. For this we inductively define $\xi_\alpha < \sigma^+$ for $\alpha < \sigma$ so that $g (\alpha, a_{\xi_\alpha}) = g (\alpha, a_\delta)$ whenever $\xi_\alpha < \delta < \sigma^+$. Thus suppose that $\xi_\beta$ has been defined for each $\beta < \alpha$. 

\medskip

{\bf Claim.}  There is $\xi < \sigma^+$ such that $g (\alpha, a_\xi) = g (\alpha, a_\delta)$ whenever $\xi < \delta < \sigma^+$. 

\smallskip

{\bf Proof of the claim.} Suppose otherwise. Inductively define $\gamma_n$ for $n < \omega$ so that
\begin{itemize}
\item $\gamma_0 = \sup \{\xi_\beta : \beta < \alpha\}$.
\item $\gamma_{n + 1} > \gamma_n$ and $g (\alpha, a_{\gamma_{n + 1}}) \not= g (\alpha, a_{\gamma_n})$.
\end{itemize}
Then $g (\alpha, a_{\gamma_0}) > g (\alpha, a_{\gamma_1}) > g (\alpha, a_{\gamma_2}) \cdots$. This contradiction completes the proof of the claim. 

\medskip

Using the claim, we let $\xi_\alpha =$ the least  $\xi < \sigma^+$ such that $g (\alpha, a_\xi) = g (\alpha, a_\delta)$ whenever $\xi < \delta < \sigma^+$. 

Finally pick $\gamma, \delta$ so that $\sup \{\xi_\eta : \eta < \sigma \} \leq \gamma < \delta < \sigma^+$. Then $g (\eta, a_\xi) = g (\eta, a_\delta)$ for all $\eta < \sigma$. Hence $F (a_\gamma, a_\delta) = 0$, which yields the desired contradiction.
 \hfill$\square$

\medskip

The remainder of the section is devoted to partitions of $\kappa \times P_\kappa (\lambda)$.

\medskip
 
\begin{Def} Given two collections $X$ and $Y$ of subsets of $P_\kappa (\lambda)$, an ideal $J$ on $P_\kappa (\lambda)$, and a cardinal $\rho$ with $0 < \rho \leq \kappa$, $X \xrightarrow[\kappa]{Y} (J^+, \rho)^2$ means that for any $F : \kappa \times P_\kappa (\lambda) \rightarrow 2$ and any $A \in X$, there is either $B \in J^+ \cap P (A)$ such that $\{ b \in B : F (\sup (a \cap \kappa), b) \not= 0\} \in Y$ for all $a \in B$, or an increasing sequence $\langle c_i : i < \rho \rangle$ in $(A, \subset)$ such that $F (\sup (c_i \cap \kappa), c_j) = 1$ whenever $i < j < \rho$.  
\end{Def}

\begin{Obs} $X \xrightarrow{J} (J^+, \rho)^2$ implies $X \xrightarrow[\kappa]{J} (J^+, \rho)^2$. 
\end{Obs}

\begin{Obs} Assuming $J$ is fine, the following are equivalent :
\begin{enumerate}[(i)]
\item $X \xrightarrow[\kappa]{J} (J^+, \rho)^2$ holds.
\item For any $G : \kappa \times P_\kappa (\lambda) \rightarrow 2$ and any $A \in X$, there is either $B \in J^+ \cap P (A)$ such that $\{ b \in B : G (\sup (a \cap \kappa), b) \not= 0\} \in J$ for all $a \in B$, or an increasing sequence $\langle c_i : i < \rho \rangle$ in $(A, \subset)$ such that $\sup (c_i \cap \kappa) < \sup (c_j \cap \kappa)$ and $G (\sup (c_i \cap \kappa), c_j) = 1$ whenever $i < j < \rho$.  
\end{enumerate}
\end{Obs}

{\bf Proof.}  \hskip0,4cm  (i) $\rightarrow$ (ii) : Given $G : \kappa \times P_\kappa (\lambda) \rightarrow 2$ and $A \in X$, define $F : \kappa \times P_\kappa (\lambda) \rightarrow 2$ by  : $F (\alpha, b) = 1$ if and only if $G (\alpha, b) = 1$ and $\alpha < \sup (a \cap \kappa)$.

\hskip0,4cm  (ii) $\rightarrow$ (i) : Trivial.
\hfill$\square$

\begin{Def} An ideal $J$ on $P_\kappa (\lambda)$ is $\kappa$\emph{-normal} if for any $A \in J^+$ and any $f : A \rightarrow \kappa$ with the property that $f (a) \in a$ for all $a \in A$, there is $B \in J^+ \cap P (A)$ such that $f$ is constant on $B$.

We let $NS^\kappa_{\kappa, \lambda}$ denote the smallest $\kappa$-normal, fine ideal on $P_\kappa (\lambda)$.
 \end{Def}

\begin{Def} We let $\Omega_{\kappa, \lambda}$ denote the set of all $a \in P_\kappa (\lambda)$ such that $a \cap \kappa$ is an infinite limit ordinal. 
\end{Def}

\begin{fact} {\rm (\cite{MPS1})} $\Omega_{\kappa, \lambda} \in (NS^\kappa_{\kappa, \lambda})^\ast$. 
\end{fact}

\begin{Def} An ideal $J$ on $P_\kappa (\lambda)$ is a \emph{weak $\pi$-point} if for any $A \in J^+$ and any $f : \kappa \rightarrow J$, there is $B \in J^+ \cap P (A)$ such that $B \cap f (\alpha) \in I_{\kappa, \lambda}$ for every $\alpha \in \kappa$.

$J$ is a \emph{weak $\chi$-point} if for any $A \in J^+$ and any $g : \kappa \rightarrow P_\kappa (\lambda)$, there is $B \in J^+ \cap P (A)$ such that $g (\sup (a \cap \kappa) \subseteq b$ for all $a, b \in B$ with $\sup (a \cap \kappa) < \sup (b \cap \kappa)$.
\end{Def}

\begin{Obs}
\begin{enumerate}[(i)]
\item  $I_{\kappa, \lambda}$ is a weak $\pi$-point. 
\item Any weak $\chi$-point is fine.
\item If $J$ is fine and $\kappa$-normal, then it is both a weak $\pi$-point and a weak $\chi$-point.
\end{enumerate}
\end{Obs}

%\begin{fact} {\rm (\cite{normal})} Suppose that $J$ is disjointly $(\sigma, \lambda^{< \kappa})$-distributive, where $J$ is a $\kappa$-normal, fine ideal on $P_\kappa (\lambda)$, and $\sigma$ is an infinite cardinal. Then $\{a \in \Omega_{\kappa, \lambda} : \cf (a \cap \kappa) \leq \sigma\} \in J$. \end{fact}

\begin{Def} Given a collection $X$ of subsets of $P_\kappa (\lambda)$, an ideal $J$ on $P_\kappa (\lambda)$, and a cardinal $\rho$ with $0 < \rho \leq \kappa$, $X \xrightarrow[\kappa]{\prec} (J^+, \rho)^2$ means that for any $F : \kappa \times P_\kappa (\lambda) \rightarrow 2$ and any $A \in X$, there is either $B \in J^+ \cap P (A)$ such that $F (\sup (a \cap \kappa), b) = 0$ whenever $a, b \in B$ are such that $a \subset b$ and $\sup (a \cap \kappa) < \sup (b \cap \kappa)$, or an increasing sequence $\langle c_i : i < \rho \rangle$ in $(A, \subset)$ such that $\sup (c_i \cap \kappa) < \sup (c_j \cap \kappa)$ and $F (\sup (c_i \cap \kappa), c_j) = 1$ whenever $i < j < \rho$.  
\end{Def}

\begin{fact} {\rm (\cite{normal})} Let $\rho$ be an uncountable cardinal less than or equal to $\kappa$, and $X = \{a \in \Omega_{\kappa, \lambda} : \cf (a \cap \kappa) < \rho \}$. Then $X \xrightarrow[\kappa]{\prec} (J^+, \rho)^2$ fails. 
\end{fact}

\begin{Obs} Suppose that $J^+ \xrightarrow[\kappa]{J} (J^+, \rho)^2$ holds, where $J$ is both a weak $\pi$-point and a weak $\chi$-point, and $\rho$ is a cardinal with $0 < \rho \leq \kappa$. Then $J^+ \xrightarrow[\kappa]{\prec} (J^+, \rho)^2$ holds. 
\end{Obs}

{\bf Proof.}  Use Observation 5.8.
\hfill$\square$

\begin{Obs} Suppose that $J^+ \xrightarrow[\kappa]{J} (J^+, \rho)^2$ holds, where $J$ is a $\kappa$-normal, fine ideal on $P_\kappa (\lambda)$, and $\rho$ is a cardinal with $\omega < \rho < \kappa$. Then $\{a \in \Omega_{\kappa, \lambda} : \cf (a \cap \kappa) < \rho \} \in J$. 
\end{Obs}

{\bf Proof.}  By Fact 5.15 and Observations 5.13 and 5.16.
\hfill$\square$

\medskip

Let us now concentrate on $J^+ \xrightarrow[\kappa]{J} (J^+, \kappa)^2$.

\medskip

\begin{Obs} Suppose that $\{P_\kappa (\lambda)\} \xrightarrow[\kappa]{J} (J^+, \kappa)^2$ holds, where $J$ is a $\kappa$-complete, fine ideal on $P_\kappa (\lambda)$. Then $\kappa$ is weakly compact. 
\end{Obs}

{\bf Proof.}  Given $f : \kappa \times \kappa \rightarrow 2$, define $F : \kappa \times P_\kappa (\lambda) \rightarrow 2$ by : $F (\alpha, b) = 1$ if and only if $\alpha <  \sup (b \cap \kappa)$ and $f (\alpha, \sup (b \cap \kappa)) = 1$. 
 
\underline {Case 1} : There is $B \in J^+ $ such that $T_a \in J$ for all $a \in B$, where $T_a = \{ b \in B : F (\sup (a \cap \kappa), b) \not= 0\}$. Proceeding by induction, define $a_\alpha \in A$ for $\alpha < \kappa$ so that for any $\beta < \alpha$, $a_\alpha \notin T_{a_\beta}$ and $\sup (a_\beta \cap \kappa) <  \sup (a_\alpha \cap \kappa)$. Then clearly, $f (\sup (a_\beta \cap \kappa),  \sup (a_\alpha \cap \kappa)) = 0$ whenever $\beta < \alpha < \kappa$.

\underline {Case 2} : There is an increasing sequence $\langle c_i : i < \kappa \rangle$ in $(A, \subset)$ such that $F (\sup (c_i \cap \kappa), c_j) = 1$ whenever $i < j < \kappa$. Then given $i < j < \kappa$, we have $\sup (c_i \cap \kappa) <  \sup (c_j \cap \kappa)$ and $f (\sup (c_i \cap \kappa), \sup (c_j \cap \kappa)) = 1$. 
\hfill$\square$
 
\begin{Def} Given an ideal $J$ on $P_\kappa (\lambda)$, we let $J \upharpoonright \kappa$ denote the set of all $X \subseteq \kappa$ such that $\{a \in P_\kappa (\lambda) : \sup (a \cap \kappa ) \in X\}$ lies in $J$.
\end{Def}

\begin{Def} An ideal $J$ on $\kappa$ is \emph{weakly selective} if for any $A \in J^+$ and any partition $Q$ of $A$ into sets in $J$, there is $B \in J^+ \cap P (A)$ such that $\vert B \cap W \vert \leq 1$ for every $W \in Q$.
\end{Def}

\begin{Obs}
\begin{enumerate}[(i)]
\item  $J \upharpoonright \kappa$ is an ideal on $\kappa$.
\item If $J$ is fine, then so is $J \upharpoonright \kappa$.
\item If $J$ is $\kappa$-complete, then so is  $J \upharpoonright \kappa$.
\item If $J$ is both a weak $\pi$-point and a weak $\chi$-point, then  $J \upharpoonright \kappa$ is weakly selective.
\item If $J$ is fine and $\kappa$-normal, then  $J \upharpoonright \kappa$ is normal.
\item $J \upharpoonright \kappa \subseteq K \upharpoonright \kappa$ for every ideal $K$ on $P_\kappa (\lambda)$ extending $J$.
\end{enumerate}
 \end{Obs}
 
 \begin{fact}  {\rm (\cite{Men74})} $NS_{\kappa, \lambda} \upharpoonright \kappa = NS_\kappa$.
\end{fact}

\begin{Def} Given an ideal $J$ on a set $X$, we let $cof (J)$ denote the least size of any ${\cal B} \subseteq J$ such that $J = \bigcup_{B \in {\cal B}} P (B)$.

Assuming that $J$ is $\kappa$-complete, but not $\kappa^+$-complete, we let $\overline{cof} (J)$ denote the least size of any ${\cal B} \subseteq J$ such that for any $A \in J$, there is $b \in P_\kappa (\cal{B})$ with $A \subseteq \bigcup b$.
 \end{Def}
 
\begin{Def} The \emph{dominating number} $\mathfrak{d}_\kappa$ denotes the least size of any $F \subseteq {}^\kappa \kappa$ with the property that for any $g \in {}^\kappa \kappa$, there is $f \in F$ such that $g (\alpha) < f (\alpha)$ for every $\alpha \in \kappa$.

We let $\overline{\mathfrak{d}}_\kappa$ denote the least size of any $F \subseteq {}^\kappa \kappa$ with the property that for any $g \in {}^\kappa \kappa$, there is $x \in P_\kappa (F)$ such that $g (\alpha) < \sup \{f (\alpha) : f \in x\}$ for every $\alpha \in \kappa$.
\end{Def}

\begin{fact}
\begin{enumerate}[(i)]
\item {\rm (\cite{Landver})} $cof (NS_\kappa) = \mathfrak{d}_\kappa$.
\item{\rm (\cite{MRS})}  $\overline{cof} (NS_\kappa) = \overline{\mathfrak{d}}_\kappa$.
\end{enumerate}
 \end{fact}

\begin{Obs}
\begin{enumerate}[(i)]
\item  Let $J$ be a $\kappa$-complete, fine ideal on $\kappa$ such that $\overline{cof} (J) \leq \lambda$. Then $J = (I_{\kappa, \lambda} \vert A) \upharpoonright \kappa$ for some $A \in I^+_{\kappa, \lambda}\cap P (\Omega_{\kappa, \lambda})$.
\item Suppose $\overline{\mathfrak{d}}_\kappa \leq \lambda$. Then $NS_\kappa = (I_{\kappa, \lambda} \vert C) \upharpoonright \kappa$ for some $C \in NS^\ast_{\kappa, \lambda}$.
\item  Let $J$ be a $\kappa$-complete, fine ideal on $\kappa$ such that $\overline{cof} (J) \leq \lambda$. Then $J = (I_{\kappa, \lambda} \vert A) \upharpoonright \kappa$ for some $A \in I^+_{\kappa, \lambda}$ such that $(A, \subset)$ and $(P_\kappa (\lambda), \subset)$ are isomorphic.
\end{enumerate}
 \end{Obs}
 
 {\bf Proof.}  (i) : Pick {${\cal B} \subseteq J$ with $\vert {\cal B} \vert = \overline{cof} (J)$ such that for any $B \in J$, there is $b \in P_\kappa ({\cal B})$ with $B \subseteq \bigcup b$. Select $v \subseteq \lambda \setminus \kappa$ with $\vert v \vert = \vert \cal{B} \vert$, and a bijection $h : v \rightarrow \cal{B}$. Let $A$ be the set of all $a \in \Omega_{\kappa, \lambda}$ such that $a \cap \kappa \notin h (\alpha)$ for all $\alpha \in v \cap a$.
 
\medskip

{\bf Claim.}  Let $X \in J^+$. Then $\{a \in A : a \cap \kappa \in X\} \in I^+_{\kappa, \lambda}$. 

\smallskip

{\bf Proof of the claim.} Fix $c \in P_\kappa (\lambda)$. Pick $\delta \in X \setminus (\bigcup_{\alpha \in v \cap c} h (\alpha))$ with $c \cap \kappa \subseteq \delta$. Set $a = \delta \cup (c \setminus \kappa)$. Then clearly, $v \cap c = v \cap a$. Moreover, $c \subseteq a$ and $a \in A$, which completes the proof of the claim. 

\medskip
 
By the claim, $A \in I^+_{\kappa, \lambda}$, and moreover $(I_{\kappa, \lambda} \vert A) \upharpoonright \kappa \subseteq J$. For the reverse inclusion, fix $X \in J$. We may find $e \in P_\kappa (v)$ such that $X \subseteq \bigcup_{\alpha \in e} h (\alpha)$. Then $a \cap \kappa \notin X$ for all $a \in A$ with $e \subseteq a$. Hence $X \in (I_{\kappa, \lambda} \vert A) \upharpoonright \kappa$.  

(ii) : By Fact 5.11, 5.22 and 5.25 and the proof of (i).

(iii) :  Let ${\cal B}$, $v$  and $h$ be as in the proof of (i), and let $A$ be the set of all $a \in P_\kappa (\lambda)$ such that 
\begin{itemize}
\item $\sup (a \cap \kappa) \in a$.
\item $\sup (a \cap \kappa) \notin h (\alpha)$ for all $\alpha \in v \cap a$.
\end{itemize}
 
\medskip

{\bf Claim.}  Let $X \in J^+$. Then $\{a \in A : \sup (a \cap \kappa) \in X\} \in I^+_{\kappa, \lambda}$. 

\smallskip

{\bf Proof of the claim.} Given $c \in P_\kappa (\lambda)$, pick $\delta \in X \setminus (\bigcup_{\alpha \in v \cap c} h (\alpha))$ with $\sup (c \cap \kappa) < \delta$. Set $a = c \cup \{\delta\}$. Then clearly, $c \subseteq a$, $a \in A$ and $\sup (a \cap \kappa) = \delta$, which completes the proof of the claim. 

\medskip
 
By the claim, $A \in I^+_{\kappa, \lambda}$, and moreover $(I_{\kappa, \lambda} \vert A) \upharpoonright \kappa \subseteq J$. For the reverse inclusion, proceed as in the proof of (i). Finally, define $\varphi : P_\kappa (\lambda) \rightarrow \kappa$ and $\psi : P_\kappa (\lambda) \rightarrow P_\kappa (\lambda)$ by $\varphi (a) =$ the least $\delta > \sup (a \cap \kappa)$ such that $\delta \notin \bigcup_{\alpha \in v \cap a} h (\alpha)$, and $\psi (a) = a \cup \{\varphi (a)\}$. It is easy to see that $\psi$ is an isomorphism from $(P_\kappa (\lambda), \subset)$ onto $(A, \subset)$.
\hfill$\square$

\medskip

Notice that if $(A, \subset)$ and $(P_\kappa (\lambda), \subset)$ are isomorphic and, say, $\{P_\kappa (\lambda)\} \xrightarrow[\kappa]{I_{\kappa, \lambda}} (I_{\kappa, \lambda}^+, \kappa)^2$ holds, then so does $\{A\} \xrightarrow[\kappa]{I_{\kappa, \lambda}} (I_{\kappa, \lambda}^+, \kappa)^2$.

\medskip

 \begin{Def} Given a collection $W$ of subsets of $\kappa$, an ideal $K$ on $\kappa$, and a cardinal $\rho$ with $0 < \rho \leq \kappa$, $W \rightarrow (K^+, \rho)^2$ means that for any $f : \kappa \times \kappa  \rightarrow 2$ and any $A \in W$, there is either $B \in K^+ \cap P (A)$ such that $f (\alpha, \beta) = 0$ for any $\alpha < \beta$ in $B$, or an increasing sequence $\langle \gamma_i : i < \rho \rangle$ in $(A, <)$ such that $f (\gamma_i, \gamma_j) = 1$ whenever $i < j < \rho$.  
\end{Def}
 
 \begin{Pro} Suppose that $J^+ \xrightarrow[\kappa]{J} (J^+, \kappa)^2$ holds, where $J$ is both a weak $\pi$-point and a weak $\chi$-point. Then  $(J \upharpoonright \kappa)^+ \rightarrow ((J \upharpoonright \kappa)^+, \kappa)^2$ holds. 
\end{Pro}

{\bf Proof.}  Let $X \in (J \upharpoonright \kappa)^+$ and $f : \kappa \times \kappa \rightarrow \kappa$ be given. Put $A = \{a \in P_\kappa (\lambda) : \sup (a \cap \kappa) \in X\}$, and define $F : \kappa \times A \rightarrow 2$ by : $F (\alpha, b) = 1$ if and only if $\alpha < \sup (b \cap \kappa)$ and $f (\alpha, \sup (b \cap \kappa)) = 1$.

\underline {Case 1} : There is $B \in J^+ \cap P (A)$, and $Z_a \in J$ for $a \in B$ such that $F (\sup (a \cap \kappa), b) = 0$ whenever $a, b \in B$ and $b \notin Z_a$. Set $T = \{\sup (a \cap \kappa) : a \in B\}$. For $\alpha \in T$, pick $a_\alpha \in B$ with $\sup (a_\alpha \cap \kappa) = \alpha$. There must be $S \in J^+ \cap P (B)$ such that for any $\alpha \in T$, and any $a, b \in S$ with $\sup (a \cap \kappa) = \alpha < \sup (b \cap \kappa)$,  $b \notin Z_{a_\alpha}$. Now put $Y = \{\sup (b \cap \kappa) : b \in S\}$. Then $Y \in (J \upharpoonright \kappa)^+ \cap P (X)$. Given $\alpha < \beta$ in $Y$, we may find $a, b \in S$ such that $\sup (a \cap \kappa) = \alpha$ and $\sup (b \cap \kappa) = \beta$. Then $b \notin Z_{a_\alpha}$. Since $a_\alpha, b \in B$, it follows that 

\centerline{$0 = F (\sup (a_\alpha \cap \kappa), b) = f (\alpha, \sup (b \cap \kappa)) = f (\alpha, \beta)$.} 
 
\underline {Case 2} : There is an increasing sequence $\langle a_\delta : \delta < \kappa \rangle$ in $(A, \subset)$ such that $F (\sup (a_\gamma \cap \kappa), a_\delta) = 1$ whenever $\gamma < \delta < \kappa$. Then clearly by definition of $F$, given $\gamma < \delta < \kappa$, we have $\sup (a_\gamma \cap \kappa) < \sup (a_\delta \cap \kappa)$, and moreover $f (\sup (a_\gamma \cap \kappa), \sup (a_\delta \cap \kappa)) = 1$.
 \hfill$\square$
 
\begin{Def} We let $NWC_\kappa$ denote the set of all $A \subseteq \kappa$ such that $\{A\} \rightarrow (NS_\kappa^+, \kappa)^2$ does not hold.
\end{Def}
 
\begin{fact} \begin{enumerate}[(i)]
\item  {\rm (\cite{Jim})} Assuming $\kappa$ is weakly compact, $NWC_\kappa$ is the smallest normal ideal $K$ on $\kappa$ such that $K^+ \rightarrow (K^+, \kappa)^2$.
\item  {\rm (\cite{Jim})} The set of all those cardinals less than $\kappa$ that are not inaccessible belongs to $NWC_\kappa$.
\item  {\rm (\cite{normal})} Let $J$ be a $\kappa$-normal, fine ideal on $P_\kappa (\lambda)$. Then $J^+ \xrightarrow[\kappa]{J} (J^+, \kappa)^2$ holds if and only if $J \upharpoonright \kappa$ extends $NWC_\kappa$.
\end{enumerate}
\end{fact}

\begin{fact} {\rm (\cite{MRS}, \cite{LCCN})}  The following are equivalent :
\begin{enumerate}[(i)]
\item $\overline{\mathfrak{d}}_\kappa \leq \lambda = \cov (\lambda, \kappa^+, \kappa^+, \kappa)$. 
\item $I_{\kappa, \lambda} \vert C$ is $\kappa$-normal for some $C \in NS_{\kappa, \lambda}^\ast$.
\end{enumerate}
\end{fact}

\medskip

Notice that by Shelah's Revised GCH theorem \cite[Conclusion 1.2]{SheRGCH}, for any uncountable strong limit cardinal $\tau$, there is $\theta < \tau$ with the property that if $\theta \leq \kappa < \tau < \lambda$, then $\cov (\lambda, \kappa^+, \kappa^+, \kappa) = \lambda$.

\medskip

\begin{Pro} Suppose $\overline{\mathfrak{d}}_\kappa \leq \lambda = \cov (\lambda, \kappa^+, \kappa^+, \kappa)$. Then the following hold.
\begin{enumerate}[(i)]
\item   Let $\sigma$ be an infinite cardinal. Then for any $D \in NS^\ast_{\kappa, \lambda}$, $J^+ \xrightarrow[\kappa]{J} (J^+, \sigma^+)^2$ does not hold, where $J = I_{\kappa, \lambda} \vert \{a \in D \cap \Omega_{\kappa, \lambda} : \cf (a \cap \kappa) \leq \sigma \}$.
\item  For any $D \in NS^\ast_{\kappa, \lambda}$ and any $X \in NS^+_\kappa \cap NWC_\kappa$, $J^+ \xrightarrow[\kappa]{J} (J^+, \kappa)^2$ does not hold, where $J = I_{\kappa, \lambda} \vert \{a \in D : \sup (a \cap \kappa) \in X\}$.
\end{enumerate}
 \end{Pro}

{\bf Proof.}  (ii) : Assume toward a contradiction that $J^+ \xrightarrow[\kappa]{J} (J^+, \kappa)^2$ holds. By Fact 5.31, there is $C \in NS^\ast_{\kappa, \lambda}$ such that $I_{\kappa, \lambda} \vert C$ is $\kappa$-normal. Put $K = J \vert (C \cap \Omega_{\kappa, \lambda})$. Then clearly, $K$ is $\kappa$-normal and $K^+ \xrightarrow[\kappa]{K} (K^+, \kappa)^2$ holds. Hence by Observation 5.21  and Proposition 5.28, $K \upharpoonright \kappa$ is a normal, fine ideal on $\kappa$ such that $(K \upharpoonright \kappa)^+ \rightarrow ((K \upharpoonright \kappa)^+, \kappa)^2$ holds. By Fact 5.30, it follows that $X \in K \upharpoonright \kappa$. Contradiction.

 (i) : Proceed as in the proof of (ii), but this time appeal to Observation 5.17. 
 \hfill$\square$

\bigskip

\section{Balanced partition properties}

\medskip

We now turn to stronger partition properties that (in the case when $\cf (\lambda) \not= \kappa$) will imply that $\cov (\lambda, \kappa^+, \kappa^+, \kappa) = \lambda$.

\begin{Def} Given $2 \leq n < \omega$, two collections $X$ and $Y$ of subsets of $P_\kappa (\lambda)$, an ideal $J$ on $P_\kappa (\lambda)$, and a cardinal $\rho$, $X \xrightarrow{Y} (J^+)^n_\rho$ means that for any $F : (P_\kappa (\lambda))^n \rightarrow \rho$ and any $A \in X$, there is $i < \rho$,  $B \in J^+ \cap P (A)$, and $Z_a \in Y$ for $a \in B$ such that $F (a_1, a_2, \cdots, a_n) = i$ whenever $a_1, a_2, \cdots, a_n \in B$ and $a_{m + 1} \notin Z_{a_1} \cup Z_{a_2} \cup \cdots \cup Z_{a_m}$ for $1 \leq m < n$. 

$X \xrightarrow[\kappa]{Y} (J^+)^2_\rho$ means that for any $F : \kappa \times P_\kappa (\lambda) \rightarrow \rho$ and any $A \in X$, there is $i < \rho$ and  $B \in J^+ \cap P (A)$ such that $\{ b \in B : F (\sup (a \cap \kappa), b) \not= i\} \in Y$ for all $a \in B$.  
\end{Def}

\begin{Obs}  \begin{enumerate} [\rm (i)]
\item Suppose that  $J^+ \xrightarrow{J} (J^+)^2_2$ holds. Then $J^+ \xrightarrow{J} (J^+)^2_n$ holds for every $n$ with $0 < n < \omega$.
\item $X \xrightarrow{J} (J^+)^2_\rho$ implies $X \xrightarrow[\kappa]{J} (J^+)^2_\rho$.
\item Suppose that $J$ is fine and $(\kappa, 2)$-distributive. Then $J^+ \xrightarrow[\kappa]{J} (J^+)^2_\rho$ whenever $0 < \rho < \kappa$.
\end{enumerate}
\end{Obs}

\begin{Obs} Suppose that $\{P_\kappa (\lambda)\} \xrightarrow{J} (J^+)^2_\rho$ holds, where $\rho$ is an infinite cardinal and $J$ is a fine ideal on $P_\kappa (\lambda)$. Then $PS^+ (\rho^+, \kappa, \nu)$ holds for any cardinal $\nu$ with $\kappa \leq \nu \leq \lambda$.
 % $\{P_\kappa (\lambda)\} \xrightarrow{I_{\kappa, \lambda}} (I_{\kappa, \lambda}^+)^2_\omega$ holds. Then $PS^+(\tau, \kappa)$ (respectively $PS(\tau, \kappa)$) assert the following: for $\beta \in \kappa$, let $Q_\beta$ be a partition of $\kappa$ with $\vert Q_\beta \vert <\tau$. Then there is a cofinal subset $B$ of $\kappa$ and $h \in \prod_{\beta \in \kappa}Q_\beta$ such that for any $\alpha,\beta\in B$, we have $\{\gamma \in h(\alpha)\cap h(\beta);\ \max \{\alpha,\beta\} \leq \gamma\} \neq \emptyset$ (respectively, there is $\zeta \in \kappa$ such that $\max \{\alpha,\beta\} \leq \zeta$ and we have $\{\gamma \in h(\alpha) \cap h(\zeta);\ \zeta \leq \gamma\} \not= \emptyset$ and $\{\delta \in h(\beta)\cap h(\zeta);\ \zeta\leq \delta\} \not= \emptyset$). 
\end{Obs}

{\bf Proof.}  Let $\nu$ be a cardinal with $\kappa \leq \nu \leq \lambda$, and for $x \in P_\kappa (\nu)$, $Q_x$ be a partition of $P_\kappa (\nu)$ with $\vert Q_x \vert \leq \rho$. For $x \in P_\kappa (\nu)$, let $\langle Q^\xi_x : \xi < \vert Q_x \vert \rangle$ be a one-to-one enumeration of $Q_x$. Define $F : P_\kappa (\lambda) \times P_\kappa (\lambda) \rightarrow \rho$ by : $F (a, c) = \xi$ if and only if $c \cap \nu \in Q^\xi_{a \cap \nu}$. There must be $\xi < \rho$ and  $B \in J^+$ such that $T_a \in J$ for all $a \in B$, where $T_a = \{ c \in B : F (a, c) \not= \xi\}$. Set $X = \{a \cap \nu : a \in B\}$. Notice that $X \in I^+_{\kappa, \nu}$. Given $x, y \in X$, select $a, b \in B$ such that $x = a \cap \nu$ and $y = b \cap \nu$. Pick $c \in B \setminus (T_a \cup T_b)$ such that $a \cup b \subseteq c$. Then obviously, $x \cup y \subseteq c \cap \nu$. Moreover, $c \cap \nu \in Q^\xi_x \cap Q^\xi_y$.  
 \hfill$\square$

\begin{Obs} Suppose that $\{P_\kappa (\lambda)\} \xrightarrow{J} (J^+)^2_\omega$ holds, where $J$ is a $\kappa$-complete, fine ideal on $P_\kappa (\lambda)$. Then the following hold :
 \begin{enumerate} [\rm (i)]
 \item Assume $\cf (\lambda) \not= \kappa$. Then $\cov (\lambda, \kappa^+, \kappa^+, \kappa) = \lambda$.
\item Let $\tau$ be a cardinal with $\kappa^+ \leq \tau < \lambda$. Then $u (\kappa^+, \tau)$ equals $\tau$ if $\cf (\tau) > \kappa$, and $\tau^+$ otherwise.
\end{enumerate}
\end{Obs}

{\bf Proof.} For any cardinal $\nu$ such that $\kappa < \nu < \lambda$ and $\cf (\nu) = \omega$, we have that $PS^+ (\omega_1, \nu^+)$ holds by Proposition 3.5 and Observation 6.3, and hence that $\cov (\nu, \nu, \omega_1, 2) = \nu^+$ by Observation 4.6. Furthermore by Observation 5.18,  $\kappa$ is inaccessible. Now for (i),  apply Observation 4.8 (i). To obtain (ii), appeal to Facts 4.4 and 4.7 if $\cf (\tau) \leq \kappa$, and to Observation 4.8 (i) otherwise. 
\hfill$\square$

\medskip

The remainder of the section is, just like the end of the previous section, devoted to negative partition relations.

\medskip

\begin{Def} Given two collections $X$ and $Y$ of subsets of $P_\kappa (\lambda)$, an ideal $J$ on $P_\kappa (\lambda)$ and two nonzero cardinals $\sigma$ and  $\rho$, the \emph{square bracket partition relation} $X \xrightarrow[\sigma]{Y} [J^+]^2_\rho$ means that for any $F : \sigma \times P_\kappa (\lambda) \rightarrow \rho$ and any $A \in X$, there is $B \in J^+ \cap P (A)$ and $\xi \in \rho$ such that $\{ b \in B : F (\sup (a \cap \sigma), b) = \xi\} \in Y$ for all $a \in B$.
\end{Def}

\begin{Def} For any cardinal $\chi$ with $\kappa \leq \chi \leq \lambda$, we let $S^\chi_{\kappa, \lambda} = \{a \in P_\kappa (\lambda) : \vert a \cap \chi \vert = \vert a \cap \kappa \vert\}$.
\end{Def}

\begin{fact} {\rm (\cite{square})}  
\begin{enumerate} [\rm (i)]
\item Suppose that $\kappa$ is weakly inaccessible and $2^\kappa \leq \lambda = \cov (\lambda, \kappa^+, \kappa^+, \kappa)$. Then $\{ C \cap S^\lambda_{\kappa, \lambda}\} \xrightarrow[\kappa]{I_{\kappa, \lambda}} [I_{\kappa, \lambda}^+]^2_\lambda$ fails for some $C \in NS^\ast_{\kappa, \lambda}$. 
\item Suppose that $\kappa$ is weakly Mahlo and $\overline{\mathfrak{d}}_\kappa \leq \lambda$, and let $A$ be the set of all $a \in S_{\kappa, \lambda} ^{\overline{\mathfrak{d}}_\kappa}$ such that $a \cap \kappa$ is a regular infinite cardinal. Then there is $D \in NS_{\kappa, \lambda}^\ast$ and $F : \kappa \times P_\kappa (\lambda) \rightarrow \lambda$ with the property that for any $B \in (NS_{\kappa, \lambda}^\kappa)^+ \cap P (A \cap D)$ and any $\xi \in \lambda$, one may find $a, b \in B$ with $a \cap \kappa < b \cap \kappa$ and $F (a \cap \kappa, b) = \xi$.
\end{enumerate}
\end{fact}

\begin{Obs}  Suppose that $\kappa$ is weakly Mahlo, $\overline{\mathfrak{d}}_\kappa \leq \lambda$, and $J$ is a $\kappa$-normal, fine ideal on $P_\kappa (\lambda)$ such that $J^+ \xrightarrow[\kappa]{J} [J^+]^2_\lambda$ holds. Then $A \cap C \in J$ for some $C \in NS_{\kappa, \lambda}^\ast$, where $A$ denotes the set of all $a \in S_{\kappa, \lambda}^{\overline{\mathfrak{d}}_\kappa}$ such that $a \cap \kappa$ is a regular infinite cardinal. 
\end{Obs}

{\bf Proof.} Suppose otherwise. Then clearly,  $A \in NS^+_{\kappa, \lambda}$. Now let $D \in NS_{\kappa, \lambda}^\ast$ and $F : \kappa \times P_\kappa (\lambda) \rightarrow \lambda$. Since $\{A \cap D\} \xrightarrow[\kappa]{J} [J^+]^2_\lambda$ and $J$ is $\kappa$-normal, there must be $B \in J^+ \cap P (A \cap D \cap \Omega_{\kappa, \lambda})$ and $\xi \in \lambda$ such that $F (a \cap \kappa, b) \not= \xi$ whenever $a, b \in B$ are such that $a \cap \kappa \in b$. This contradicts Fact 6.7.
\hfill$\square$

\begin{Cor} Suppose that $\overline{\mathfrak{d}}_\kappa \leq \lambda$ and $J^+ \xrightarrow[\kappa]{J} (J^+)^2_2$ holds, where $J$ is a $\kappa$-normal, fine ideal on $P_\kappa (\lambda)$. Then $S_{\kappa, \lambda}^{\overline{\mathfrak{d}}_\kappa} \cap C \in J$ for some $C \in NS_{\kappa, \lambda}^\ast$.
\end{Cor}

{\bf Proof.} By Observation 5.18 and Fact 5.30, $\kappa$ is weakly compact, and moreover the set of all $a \in P_\kappa (\lambda)$ such that $a \cap \kappa$ is an inaccessible cardinal lies in $J^\ast$. 
\hfill$\square$

\begin{Cor} Suppose that $\kappa$ is weakly Mahlo and $\overline{\mathfrak{d}}_\kappa \leq \lambda = \cov (\lambda, \kappa^+, \kappa^+, \kappa)$. Then for any $D \in NS^\ast_{\kappa, \lambda}$, $J^+ \xrightarrow[\kappa]{J} [J^+]^2_\lambda$ does not hold, where $J = I_{\kappa, \lambda} \vert (D \cap S_{\kappa, \lambda}^{\overline{\mathfrak{d}}_\kappa})$.
 \end{Cor}

{\bf Proof.} Assume toward a contradiction that there exists $D \in NS^\ast_{\kappa, \lambda}$ such that $J^+ \xrightarrow[\kappa]{J} [J^+]^2_\lambda$ holds, where $J = I_{\kappa, \lambda} \vert (D \cap S_{\kappa, \lambda}^{\overline{\mathfrak{d}}_\kappa})$. By Fact 5.31, we may find $C \in NS^\ast_{\kappa, \lambda}$ such that $I_{\kappa, \lambda} \vert C$ is $\kappa$-normal. Put $K = J \vert C$. Then clearly, $K$ is $\kappa$-normal, and moreover $K^+ \xrightarrow[\kappa]{K} (K^+)^2_2$ holds. This contradicts Observation 6.8.
 \hfill$\square$
 
\medskip
 
A similar result will be obtained as a variant of Proposition 5.32 (ii).
 
\medskip

\begin{Pro} Suppose that $J^+ \xrightarrow[\kappa]{J} (J^+)^2_2$ holds, where $J$ is a $\kappa$-normal, fine ideal on $P_\kappa (\lambda)$. Then $(J \upharpoonright \kappa)^+ \rightarrow ((J \upharpoonright \kappa)^+)^2$ holds. 
\end{Pro}

{\bf Proof.} By the proof of Proposition 5.28.
\hfill$\square$

 \begin{Obs} Given a $\kappa$-complete, fine ideal $J$ on $\kappa$, the following are equivalent:
 \begin{enumerate} [\rm (i)]
 \item $J$ is $(\kappa, 2)$-distributive.
 \item Let $0 < \eta \leq \kappa$, and $Q_\alpha \in MAD_d (J)$ for $\alpha < \eta$ be such that $Q_\beta \subseteq \bigcup_{W \in Q_\alpha} P (W)$ whenever $\alpha < \beta < \eta$. Then there is $B \in J^+ \cap P (A)$ and $h \in \prod_{\alpha < \eta} Q_\alpha$ such that $B \setminus h (\alpha) \in J$ for all $\alpha < \kappa$.
 \end{enumerate}
\end{Obs} 

{\bf Proof.}  \hskip0,4cm  (i) $\rightarrow$ (ii) : By Observation 5.4.

\hskip0,4cm  (ii) $\rightarrow$ (iii) : Suppose that (ii) holds. We claim that $\kappa$ is inaccessible. Suppose otherwise, and let $\nu$ be the least cardinal such that $2^\nu \geq \kappa$. 
Let $\langle X_\xi : \xi < \kappa \rangle$ be a sequence of pairwise distinct subsets of $\nu$. For $\delta < \nu$, put $A^0_\delta= \{\xi < \kappa : \delta \in X_\xi\}$ and $A^1_\delta = \kappa \setminus A^0_\delta$. Now let $Q_0 = \{\kappa\}$, and for $0 < \alpha < \nu$, $Q_\alpha = \{ \bigcap_{\delta < \alpha} A^{k (\delta)}_\delta : k \in  {}^\alpha 2\} \cap J^+$. Then $\vert \bigcap_{\alpha < \nu} h (\alpha) \vert \leq 1$ for all $h \in \prod_{\alpha < \nu} Q_\alpha$, which yields the desired contradiction. 

Now suppose that $A \in J^+$ and for $\alpha < \kappa$, $W_\alpha \in  MAD_d (J)$ with $\vert W_\alpha \vert \leq 2$. For $\alpha < \kappa$, let

\centerline{$T_\alpha = \{\bigcap_{\beta \leq \alpha} g (\beta) : g \in \prod_{\beta \leq \alpha} W_\beta\} \cap J^+$.} 

There must be $B \in J^+ \cap P (A)$ and $f \in \prod_{\alpha < \kappa} T_\alpha$ such that $B \setminus f (\alpha) \in J$ for all $\alpha < \kappa$. For $\alpha < \kappa$, let $g_\alpha \in \prod_{\beta \leq \alpha} W_\beta$ be such that $f (\alpha) = \bigcap_{\beta \leq \alpha} g_\alpha (\beta)$. Then it is easy to see that $\bigcup_{\alpha < \kappa} g_\alpha \in \prod_{\beta \leq \kappa} W_\beta$. Moreover, $B \setminus (\bigcup_{\alpha < \kappa} g_\alpha) (\beta) \in J$ for all $\beta < \kappa$.
 \hfill$\square$

\begin{Def}  Given an ideal $J$ on $\kappa$,  $2 \leq n < \omega$, a collection $X$ of subsets of $\kappa$ and an ordinal $\rho$, $X \rightarrow (J^+)^n_\rho$ means that for any $F : \kappa^n \rightarrow \rho$ and any $A \in X$, there is $i < \rho$ and  $B \in J^+ \cap P (A)$ such that $F (\alpha_1, \alpha_2, \cdots, \alpha_n) = i$ whenever $\alpha_1 <  a_2 < \cdots < \alpha_n$ are in $B$. 
\end{Def}

\begin{fact} Given a $\kappa$-complete, fine ideal $J$ on $\kappa$, the following are equivalent : 
\begin{enumerate} [\rm (i)]
\item $J^+ \rightarrow (J^+)^2_2$.
\item$J^+ \rightarrow (J^+)^n_\rho$ whenever $0 < n < \omega$ and $0 < \rho < \kappa$.
\item $J$ is $(\kappa, 2)$-distributive and weakly selective.
 \end{enumerate}
\end{fact} 

{\bf Proof.}  By Theorem 9 in \cite{Johnson} and Observation 6.12.
\hfill$\square$

%\begin{fact} {\rm(\cite{Fodor})} Any normal, fine ideal $J$ on $\kappa$ is weakly selective.\end{fact} 

\begin{Def} $\kappa$ is \emph{completely ineffable} if there exists a normal, $(\kappa, 2)$-distributive, fine ideal on $\kappa$.
\end{Def}

\begin{fact} Suppose that $\kappa$ is completely ineffable. Then there exists a smallest normal, $(\kappa, 2)$-distributive, fine ideal on $\kappa$. 
\end{fact}

{\bf Proof.} By the proof of Corollary 3 in \cite{Johnson} and Observation 6.12.
\hfill$\square$

\begin{Def} We let $NCI_\kappa$ denote the smallest normal, $(\kappa, 2)$-distributive, fine ideal on $\kappa$ if $\kappa$ is completely ineffable, and $P (\kappa)$ otherwise.
\end{Def}

\begin{Pro} Suppose $\overline{\mathfrak{d}}_\kappa \leq \lambda = \cov (\lambda, \kappa^+, \kappa^+, \kappa)$. Then for any $D \in NS^\ast_{\kappa, \lambda}$ and any $X \in NS^+_\kappa \cap NCI_\kappa$, $J^+ \xrightarrow[\kappa]{J} (J^+)^2_2$ does not hold, where $J = I_{\kappa, \lambda} \vert \{a \in D : \sup (a \cap \kappa) \in X\}$.
 \end{Pro}

{\bf Proof.} Assume toward a contradiction that there are $D \in NS^\ast_{\kappa, \lambda}$ and $X \in NS^+_\kappa \cap NCI_\kappa$ such that $J^+ \xrightarrow[\kappa]{J} (J^+)^2_2$ holds, where $J = I_{\kappa, \lambda} \vert \{a \in D : \sup (a \cap \kappa) \in X\}$. By Fact 5.31, there is $C \in NS^\ast_{\kappa, \lambda}$ such that $I_{\kappa, \lambda} \vert C$ is $\kappa$-normal. Put $K = J \vert C$. Then clearly, $K$ is $\kappa$-normal, and moreover $K^+ \xrightarrow[\kappa]{K} (K^+)^2_2$ holds. Hence by Observation 5.21, Proposition 6.11 and Fact 6.14, $K \upharpoonright \kappa$ is a normal, $(\kappa, 2)$-distributive, fine ideal on $\kappa$. It follows that $X \in J \upharpoonright \kappa$. Contradiction.
 \hfill$\square$

\begin{Pro} Suppose that $\overline{\mathfrak{d}}_\kappa \leq \lambda$ and $\cf \lambda) \not= \kappa$. Then the following hold :
\begin{enumerate}[(i)]
\item For any $D \in NS^\ast_{\kappa, \lambda}$, $J^+ \xrightarrow{J} (J^+)^2_\omega$ does not hold, where  $J = I_{\kappa, \lambda} \vert (D \cap S_{\kappa, \lambda}^{\overline{\mathfrak{d}}_\kappa})$.
\item For any $D \in NS^\ast_{\kappa, \lambda}$ and any $X \in NS^+_\kappa \cap NCI_\kappa$, $J^+ \xrightarrow{J} (J^+)^2_\omega$ does not hold, where $J = I_{\kappa, \lambda} \vert \{a \in D : \sup (a \cap \kappa) \in X\}$.
\end{enumerate}
\end{Pro}

{\bf Proof.}  By Observation 6.4, Corollary 6.10 and Proposition 6.18.
 \hfill$\square$
 
 \medskip
 
 In contrast to this, by a result of Usuba (see the proof of Theorem 1.9 in \cite{Usuba2}), it is consistent relative to a large cardinal that \say{$\kappa$ is not subtle, but $I^+_{\kappa, \lambda} \rightarrow (I^+_{\kappa, \lambda})^n_\eta$ holds for any $n < \omega$ and any $\eta < \kappa$}.

\bigskip

\section{Mild ineffability}

\medskip

\begin{Def} $\kappa$ is \emph{mildly $\lambda$-ineffable} if, given $s_a\subseteq a$ for $a \in P_\kappa (\lambda)$, there exists $S \subseteq \lambda$ with the property that for any $ b \in P_\kappa (\lambda)$, there is $a \in P_\kappa (\lambda)$ such that $b \subseteq a$ and $S \cap b = s_a \cap b$.
\end{Def}

\medskip

We will establish that if $\kappa$ is mildly $\lambda$-ineffable and $\cf (\lambda) \not= \kappa$, then $\cov (\lambda, \kappa^+, \kappa^+, \kappa) = \lambda$.  We need some preparation.

\medskip

\begin{fact} {\rm(\cite{Carr})}
\begin{enumerate}[(i)]
 \item  If $\kappa$ is mildly $\lambda$-ineffable, then it is mildly $\nu$-ineffable for any cardinal $\nu$ with $\kappa \leq \nu\leq \lambda$.
 \item $\kappa$ is mildly $\kappa$-ineffable if and only if it is weakly compact.
 \end{enumerate}
 \end{fact}
 
\begin{fact} {\rm(\cite{Laura})} Suppose that $\kappa$ is inaccessible. Then the following hold :
\begin{enumerate}[(i)]
\item $\kappa$ is mildly $\lambda$-ineffable iff $TP(\kappa,\lambda)$ holds iff  $TP^{-}(\kappa,\lambda)$ holds.
\item  $\kappa$ is mildly $\lambda^{< \kappa}$-ineffable iff  $PS^+ (\kappa, \kappa, \lambda)$ holds iff  $PS^{\ast} (\kappa, \kappa, \lambda)$ holds iff  $PS (\kappa, \kappa, \lambda)$ holds.
 \end{enumerate}
 \end{fact}
 
 \begin{fact} {\rm (\cite{Carlo})} The following are equivalent :
\begin{enumerate}[(i)]
\item $\kappa$ is mildly $\lambda$-ineffable.
\item Given $W_\alpha \subseteq P_\kappa (\lambda)$ for $\alpha < \lambda$, there is $h \in \prod_{\alpha < \lambda} \{W_\alpha, P_\kappa (\lambda) \setminus W_\alpha\}$ such that $\bigcap_{\alpha \in e} h( \alpha )\in I_{\kappa,\lambda}^+$ for every nonempty $e \in P_\kappa (\lambda)$.
\end{enumerate}
 \end{fact}

\begin{Obs} Suppose that $\kappa$ is mildly $\lambda$-ineffable, and for each $\alpha < \lambda$, let $Q_\alpha$ be a partition of $P_\kappa (\lambda)$ into less than $\kappa$ many pieces. Then there is $h \in \prod_{\alpha < \lambda} Q_\alpha$ such that $\bigcap_{\alpha \in e} h( \alpha )\in I_{\kappa,\lambda}^+$ for every nonempty $e \in P_\kappa (\lambda)$.
 \end{Obs}

{\bf Proof.}   By Fact 7.4, we may find $h \in \prod_{W \in \bigcup_{\alpha < \lambda} Q_\alpha} \{W, P_\kappa (\lambda) \setminus W\}$ such that $\bigcap_{W \in x} h(W) \in I_{\kappa,\lambda}^+$ for every nonempty $x \in P_\kappa (\bigcup_{\alpha < \lambda} Q_\alpha)$. Now given $\alpha < \lambda$, we have $\bigcap_{W \in Q_\alpha} (P_\kappa (\lambda) \setminus W) = \emptyset$, and consequently $Q_\alpha \cap ran (h) \not= \emptyset$.
\hfill$\square$
 
\medskip

$TP (\kappa, \lambda)$ may be reformulated in the same way.

\medskip

\begin{Pro} The following are equivalent :
\begin{enumerate}[(i)]
\item $TP (\kappa, \lambda)$.
\item For each $\alpha < \lambda$, let $Q_\alpha$ be a partition of $\{x \in P_\kappa (\lambda) : \alpha \in x\}$ into less than $\kappa$ many pieces. Suppose that 

\centerline{$\vert \{ \bigcap_{\alpha \in d} g (\alpha) : g \in \prod_{\alpha \in d} Q_\alpha \} \vert < \kappa$}

for any nonempty $d \in P_\kappa (\lambda)$. Then there is $h \in \prod_{\alpha < \lambda} Q_\alpha$ such that $\bigcap_{\alpha \in e} h( \alpha )\in I_{\kappa,\lambda}^+$ for every nonempty $e \in P_\kappa (\lambda)$.
\item For each $\alpha < \lambda$, let $Q_\alpha$ be a partition of $\{x \in P_\kappa (\lambda) : \alpha \in x\}$ with $\vert Q_\alpha \vert \leq 2$. Suppose that 

\centerline{$\vert \{ \bigcap_{\alpha \in d} g (\alpha) : g \in \prod_{\alpha \in d} Q_\alpha\} \vert < \kappa$}

for any nonempty $d \in P_\kappa (\lambda)$. Then there is $h \in \prod_{\alpha < \lambda} Q_\alpha$ such that $\bigcap_{\alpha \in e} h( \alpha )\in I_{\kappa,\lambda}^+$ for every nonempty $e \in P_\kappa (\lambda)$.
\end{enumerate}
 \end{Pro}

{\bf Proof.}    \hskip0,4cm  (i) $\rightarrow$ (ii) :  Assume that $TP (\kappa, \lambda)$ holds. Let $Q_\alpha$ be a partition of $\{x \in P_\kappa (\lambda) : \alpha \in x\}$ into less than $\kappa$ many pieces for $\alpha < \kappa$ such that 

\centerline{$\vert \{ \bigcap_{\alpha \in d} g (\alpha) : g \in \prod_{\alpha \in d} Q_\alpha \} \vert < \kappa$}

for any nonempty $d \in P_\kappa (\lambda)$. For $\alpha < \lambda$, let $\langle Q^\xi_\alpha : \xi < \vert Q_\alpha \vert \rangle$ be a one-to-one enumeration of $Q_\alpha$. Select a bijection $f : \kappa \times \lambda \rightarrow \lambda$. For $a \in P_\kappa (\lambda)$, define $t_a : a \rightarrow 2$ as follows. Given $\xi < \kappa$ and $\alpha < \lambda$ such that $f (\xi, \alpha) \in a$, we let $t_a (j (\xi, \alpha)) = 1$ just in case $\alpha \in a$ and $a \in Q^\xi_\alpha$. For $c \in P_\kappa (\lambda)$, let $A_c$ denote the collection of all $\alpha < \lambda$ such that $f (\xi, \alpha) \in c$ for some $\xi < \kappa$. Let $C$ be the set of all $c \in P_\kappa (\lambda)$ such that $A_c \subseteq c$. Note that $C \in NS^\ast_{\kappa, \lambda}$. 

\medskip

{\bf Claim 1.}  Let $c \in C$. Then $\vert \{t_a \vert c : c \subseteq a\} \vert <\kappa$.

\smallskip

{\bf Proof of Claim 1.} Suppose otherwise, and let $a_i \in P_\kappa (\lambda)$ for $i < \kappa$ be such that
\begin{itemize}
\item $c \subseteq a_i$ for all $i < \kappa$.
\item $t_{a_i} \vert c \not= t_{a_j} \vert c$ whenever $i < j < \kappa$.
\end{itemize}
For $i < \kappa$, define $k_i : A_c \rightarrow \kappa$ so that $a_i \in Q_\alpha^{k_i (\alpha)}$ for all $\alpha \in A_c$. Now given $i < j < \kappa$, we may find $\alpha \in A_c$ and $\xi < \kappa$ such that $f (\xi, \alpha) \in c$ and $t_{a_i} (f (\xi, \alpha)) \not= t_{a_i} (f (\xi, \alpha))$. Then it is easy to see that $k_i (\alpha) \not= k_j (\alpha)$. Hence, $k_i \not= k_j$. This contradiction completes the proof of the claim. 

\medskip

By Claim 1, we may find $T : \lambda \rightarrow 2$ such that for any $v \in P_\kappa (\lambda)$, there is $a \in P_\kappa (\lambda)$ with $v \subseteq a$ and $T \vert v = t_a \vert v$.

\medskip

{\bf Claim 2.}  Let $\alpha < \lambda$. Then $T (f (\xi, \alpha)) = 1$ for some $\xi < \vert Q_\alpha \vert$.

\smallskip

{\bf Proof of Claim 2.} Suppose otherwise. Pick $v \in P_\kappa (\lambda)$ such that $\{ \alpha \} \cup \{ f (\xi, \alpha) : \xi < \vert Q_\alpha \vert \} \subseteq v$. There must be $a \in P_\kappa (\lambda)$ such that $v \subseteq a$ and $T \vert v = t_a \vert v$. Then $a \notin Q^\xi_\alpha$ for all $\xi < \vert Q_\alpha \vert$. This contradiction completes the proof of the claim. 

\medskip

{\bf Claim 3.}  Let $\alpha < \lambda$. Then $\vert \{\xi < \lambda : T (f (\xi, \alpha)) = 1\} \leq 1$.

\smallskip

{\bf Proof of Claim 3.} Let $\xi_1, \xi_2 < \lambda$ be such that $T (f (\xi_1, \alpha)) = T (f (\xi_2, \alpha)) = 1$. Pick $v \in P_\kappa (\lambda)$ such that $\{ \alpha , f (\xi_1, \alpha), f (\xi_2, \alpha) \} \subseteq v$. There must be $a \in P_\kappa (\lambda)$ such that $v \subseteq a$ and $T \vert v = t_a \vert v$. Then $a  \in Q^{\xi_1}_\alpha \cap Q^{\xi_2}_\alpha$. Hence $\xi_1 = \xi_2$, which completes the proof of the claim. 

\medskip

Using Claims 2 and 3, define $H \in \prod_{\alpha < \lambda} \vert Q_\alpha \vert$ by $H (\alpha) =$ the unique $\xi < \vert Q_\alpha \vert$ such that $T (f (\xi, \alpha)) = 1$. Now given $e, w \in P_\kappa (\lambda) \setminus \{ \emptyset \}$, set $v = e \cup w \cup \{ f (H (\alpha), \alpha) : \alpha \in e\}$. We may find $a \in P_\kappa (\lambda)$ such that $v \subseteq a$ and $T \vert v = t_a \vert v$. Then clearly,
\begin{itemize}
\item $w \subseteq a$.
\item $a \in \bigcap_{\alpha \in e} Q^{H (\alpha)}_\alpha$.
\end{itemize} 

\medskip

\hskip0,2cm  (ii) $\rightarrow$ (iii) :  Trivial.

\medskip
\hskip0,2cm  (iii) $\rightarrow$ (i) :  Assume that (iii) holds. Let $t_a : a \rightarrow 2$ for $a \in P_\kappa (\lambda)$ be such that $\vert \{t_a \vert d : d \subseteq a\} \vert < \kappa$ for all $d \in P_\kappa (\lambda)$. For $\alpha < \lambda$ and $i < 2$, let $Q^{i}_\alpha$ be the set of all $a \in P_\kappa (\lambda)$ such that $\alpha \in a$ and $t_a (\alpha) = i$.  

\medskip

{\bf Claim.}  Let $d \in P_\kappa (\lambda) \setminus \{\emptyset\}$. Then 

\centerline{$\vert \{ \bigcap_{\alpha \in d} u (\alpha) : u \in \prod_{\alpha \in d} \{Q^0_\alpha, Q^1_\alpha\} \} \vert < \kappa$}.

\smallskip

{\bf Proof of the claim.} Suppose otherwise. Pick $g_\xi : d \rightarrow 2$ for $\xi < \kappa$ so that $\bigcap_{\alpha \in d} Q^{g_\eta (\alpha)}_\alpha \not= \bigcap_{\alpha \in d} Q^{g_\xi (\alpha)}_\alpha$ whenever $\eta < \xi < \kappa$. For $\xi < \kappa$ with $\bigcap_{\alpha \in d} Q^{g_\xi (\alpha)}_\alpha \not= \emptyset$, pick $a_\xi \in \bigcap_{\alpha \in d} Q^{g_\xi (\alpha)}_\alpha$. Notice that $t_{a_\xi} \vert d = g_\xi$. Thus $\vert \{t_{a_\xi} \vert d : \xi < \kappa\} \vert = \kappa$. This contradiction completes the proof of the claim. 

\medskip

By the claim, we may find $T : \lambda \rightarrow 2$ such that  $\bigcap_{\alpha \in e} Q_\alpha^{T (\alpha)} \in I_{\kappa,\lambda}^+$ for every nonempty $e \in P_\kappa (\lambda)$. It remains to observe that for any $e \in P_\kappa (\lambda) \setminus \{\emptyset\}$ and any $a \in \bigcap_{\alpha \in e} Q_\alpha^{T (\alpha)}$, we have $t_a \vert e = T \vert e$.

%\hskip0,2cm  (ii) $\rightarrow$ (i) :  Suppose that (ii) holds, and let $A \in J^+$, and $W_a \subseteq P_\kappa (\lambda)$ for $a \in P_\kappa (\lambda)$. For $a \in P_\kappa (\lambd It remains to observe that fa)$, set $W^{0}_a = W_a$ and $W^{1}_a = P_\kappa (\lambda) \setminus W_a$.  It is readily seen that $J^+ \xrightarrow{J} (J^+)^2$ (and hence $J^+ \xrightarrow{J} (J^+, \kappa)^2$) holds. By Observation, it follows that $\kappa$ is weakly compact (and therefore inaccessible).  Hence we may find $B \in J^+ \cap P (A)$ and $h_a : P (a) \rightarrow 2$ for $a \in B$ such that \centerline{$\{ b \in B : \exists c \subseteq a (b \notin W_c^{h_a (c)})\} \in J$} for all $a \in B$. Put $h = \bigcup_{a \in B} h_a$. Then, as is easily seen, $h$ is a func. tion from $P_\kappa (\lambda)$ to $2$. Moreover, $B \setminus W_c^{h (c)} \in J$ for any $c \in P_\kappa (\lambda)$.
\hfill$\square$

%\begin{fact} {\rm (\cite{LCCN})} Suppose that $\cf (\lambda) = \omega$ and $\kappa$ is mildly $\lambda^+$-ineffable. Then  $\cov (\lambda, \lambda, \omega_1, 2) = \lambda^+$. \end{fact}{\bf Proof.}   Suppose otherwise. Then by Fact, we may find $y_\alpha \in P_{\omega_1} (\lambda)$ for $\alpha < \lambda^+$ such that for any nonzero $\beta < \lambda^+$, there is a one-to-one $h_\beta \in \prod_{\alpha < \beta} y_\alpha$. For $\alpha < \lambda^+$, let $\langle y^n_\alpha : n < \vert y_\alpha \vert \rangle$ be a one-to-one enumeration of $y_\alpha$. Pick a bijection $j : \lambda^+ \times \omega \rightarrow \lambda^+$. For $a \in P_\kappa (\lambda^+)$, let \centerline{$s_a = a \cap \{j (\alpha, n) : \alpha < \sup a$ and $h_{\sup a} (\alpha) = y^n_\alpha\}$.}Define $\psi : P_3 (\lambda^+) \setminus \{\emptyset\} \rightarrow P_{\omega_1} (\lambda^+)$ by $\psi (e) = e \cup \{(\max e) + 1\} \cup \{j (\alpha, n) : \alpha \in e$. There must be $F \subseteq \lambda^+ \times \omega$ with the property that for any $ e \in P_3 (\lambda) \setminus \{\emptyset\}$, there is $a \in P_\kappa (\lambda^+)$ such that $\psi (e) \subseteq a$ and $j`` F \cap \psi (e) = s_a \cap \psi (e)$. Then it is easy to see that $F$ is a function with domain $\lambda^+$. Moreover, the function $g : \lambda^+ \rightarrow \lambda$ defined by $g (\alpha) = y_\alpha^{F (\alpha)}$ is one-to-one. Contradiction.\hfill$\square$
 
\begin{Def} Given a set $P$ and a $\kappa$-complete ideal $J$ on $P$, we denote by $IE_\kappa^2(J)$ the following statement : Suppose that for each $p \in P$, there is a partition $Q_p$ of $P$ with $\vert Q_p \vert < \kappa$. Then there is $h \in \prod_{p \in P} Q_p$ and a $\kappa$-complete ideal $K$ on $P$ extending $J$ such that $\ran(h) \subseteq K^\ast$.
\end{Def}

\begin{Obs} Suppose that $\kappa$ is mildly $\lambda$-ineffable and $\lambda$ is regular. Then $IE_\kappa^2(I_\lambda)$ holds. 
 \end{Obs}

{\bf Proof.}   For each $\alpha < \lambda$, let $W_\alpha$ be a partition of $\lambda$ with $\vert W_\alpha \vert < \kappa$. For $\alpha < \lambda$, put

\centerline{$Q_\alpha = \{\{a \in P_\kappa (\lambda) : \sup a \in T \} : T \in W_\alpha\}$.}

By Observation 7.5, we may find $h \in \prod_{\alpha < \lambda} W_\alpha$ such that 

\centerline{$\{a \in P_\kappa (\lambda) : \sup a \in \bigcap_{\alpha \in e} h (\alpha)\} \in I^+_{\kappa, \lambda}$}

for every nonempty $e \in P_\kappa (\lambda)$. It is easy to see that $\bigcap_{\alpha \in e} h (\alpha) \in I^+_\lambda$ for all $e \in P_\kappa (\lambda) \setminus \{\emptyset\}$. 
\hfill$\square$

\medskip

Note that if $\lambda$ is weakly compact, then by Fact 7.2 and Observation 7.8, $IE_\kappa^2(I_\lambda)$ (in fact $IE_\lambda^2(I_\lambda)$) holds. Thus the converse of Observation 7.8 does not hold.

\medskip

\begin{Obs} Suppose that $\kappa$ is mildly $\lambda$-ineffable. Then $\cov (\nu, \nu, (\cf (\nu)^+, 2) = \nu^+$ for each singular cardinal $\nu$ with $\kappa < \nu < \lambda$. 
 \end{Obs}

{\bf Proof.}  Given a singular cardinal $\nu$ with $\kappa < \nu < \lambda$, $\kappa$ is mildly $\nu^+$-ineffable by Fact 7.2, so $PS^+ ((\cf (\nu)^+, \nu^+)$ holds by Observation 7.7, and therefore $\cov (\nu, \nu, (\cf (\nu))^+, 2) = \nu^+$ by Observation 4.6. 
\hfill$\square$

\begin{fact}  {\rm(\cite{Usuba})} Suppose that  $cf(\lambda)\geq \kappa$ and $\kappa$ is mildly $\lambda$-ineffable. Then $\lambda^{<\kappa} = \lambda$.
\end{fact}

{\bf Proof.}  Since $\kappa$ is inaccessible by Fact 7.2, it follows from Fact 4.2 and Observations 4.8 (i) and 7.9 that $\lambda^{<\kappa} = u (\kappa, \lambda) = \lambda$. 
\hfill$\square$

\medskip

Usuba \cite{Usuba} asked whether $\lambda^{< \kappa} = \lambda^+$ whenever $\kappa$ is mildly $\lambda$-ineffable and $\cf (\lambda) < \kappa$. The following provides a partial answer to this question.

\medskip

\begin{Pro} \begin{enumerate}[(i)]
\item Suppose that $\omega < cf(\lambda) < \kappa$, and $\kappa$ is mildly $\nu$-ineffable for every cardinal $\nu$ with $\kappa \leq \nu < \lambda$. Then $\lambda^{<\kappa} = \lambda^+$.
\item Suppose that $\omega = cf(\lambda)$, and $\kappa$ is mildly $\nu$-ineffable for every cardinal $\nu$ with $\kappa \leq \nu < \lambda$. Then $\lambda^{<\kappa} = \cov (\lambda, \lambda, \omega_1, 2)$.
\end{enumerate}
\end{Pro}
 
{\bf Proof.}  (i) : $\kappa$ is inaccessible by Fact 7.2, so by Facts 4.2 and 4.7 and Observation 7.9, $\lambda^+ \leq \lambda^{<\kappa} = u (\kappa, \lambda) \leq \lambda^+$. 

(ii) : Use Observation 4.10.
\hfill$\square$

\begin{Pro} Suppose that $\kappa$ is mildly $\lambda$-ineffable and $\cf (\lambda) \not= \kappa$. Then $\cov (\lambda, \kappa^+, \kappa^+, \kappa) = \lambda$. 
 \end{Pro}

{\bf Proof.}   Since $\kappa$ is inaccessible by Fact 7.2, the result follows from Observations 7.9 and 4.8 (use (i) if $\cf (\lambda) \not= \omega$, and (ii) otherwise).
\hfill$\square$

\begin{fact}  {\rm(\cite{Laura})} The following are equivalent:
\begin{enumerate}[(i)]
\item $\kappa$ is mildly $\lambda^{<\kappa}$-ineffable.
\item For any set $P$ of size $\lambda^{<\kappa}$ and any $\kappa$-complete ideal $J$ on $P$, $IE_\kappa^2(J)$ holds.
\end{enumerate}
\end{fact}

\medskip

Suppose that $\kappa$ is the successor of a singular limit of $\lambda$-compact cardinals. Then by a result of Magidor and Shelah \cite{MS}, $\kappa$ has the tree property, and in fact, as shown in \cite{Laura1}, $TP(\kappa,\lambda')$ holds for every cardinal $\lambda'\geq \kappa$. We modify the proof so as to obtain the following which improves a result of \cite{Laura}.

\medskip

\begin{Pro} Suppose that $\kappa = \nu^+$, where $\nu$ is a singular limit of mildly $\lambda^{< \kappa}$-ineffable cardinals. Then $PS (\kappa, \lambda)$ holds.
\end{Pro}

{\bf Proof.} Set $\sigma= cf(\nu)$, and select an increasing sequence $\langle \nu_i : i < \sigma \rangle$ of mildly $\lambda^{< \kappa}$-ineffable cardinals with $\sigma < \nu_0$ and $\sup \{\nu_i : i < \sigma\} = \nu$. Suppose that for each $a \in P_\kappa (\lambda)$, $Q_a$ is a partition of $P_\kappa (\lambda)$ with $\vert Q_a \vert \leq \nu$. For $a \in P_\kappa (\lambda)$, pick an onto function $\psi_a$ from $Q_a$ to $\nu$, and let $W^p_a$ denote the set of all $c \in P_\kappa (\lambda)$ such that $p =$ the least $r$ such that $c \in \bigcup \psi_a``r$. By Fact 7.13, we may find $t : P_\kappa (\lambda) \rightarrow \sigma$ and a $\nu_0$-complete ideal $K$ on $P_\kappa (\lambda)$ extending $I_{\kappa, \lambda}$ such that $\{ W_a^{t (a)} : a \in P_\kappa (\lambda)\} \subseteq K^\ast$. There must be $S \in K^+$ and $p < \sigma$ such that $t$ takes the constant value $p$ on $S$. Thus $W^p_a \cap W^p_x \in I^+_{\kappa, \lambda}$ for all $a, x \in S$. Define $g: S \times S \rightarrow P_{\kappa} (\lambda)$ so that  $a \cup x \subseteq g (a, x)$ and $g (a, x) \in W^p_a \cap W^p_x$. Further define $h: S \times S \rightarrow \nu_p \times \nu_p$ so that $g (a, x) \in \psi_a (\alpha) \cap \psi_x (\beta)$, where $h (a, x) = (\alpha, \beta)$. For $a \in S$ and $(\alpha, \beta) \in \nu_p \times \nu_p$, put $T^{(\alpha, \beta)}_a = \{x \in S : h (a, x) = (\alpha, \beta)\}$.  By Fact 7.13, we may find $u : S \rightarrow \nu_p \times \nu_p$ and a $\nu_{p + 1}$-complete ideal $G$ on $P_\kappa (\lambda)$ extending $I_{\kappa, \lambda} \vert S$ such that $\{ T_a^{u (a)} : a \in S\} \subseteq G^\ast$. There must be $A \in G^+$ and $(\alpha, \beta) \in \nu_p \times \nu_p$ such that $u$ is constantly $(\alpha, \beta)$ on $A$. Pick $X \in I^+_{\kappa, \lambda} \cap P (A)$ so that $A \setminus X \in I^+_{\kappa, \lambda}$. Set $B = X$ if $A \setminus X \in G^+$, and $B = A \setminus X$ otherwise. Now given $a, b \in B$, pick $x \in T^{(\alpha, \beta)}_a \cap T^{(\alpha, \beta)}_b \cap (A \setminus B)$ with $a \cup b \subseteq x$. Then clearly, $g (a, x) \in \psi_a (\alpha) \cap \psi_x (\beta)$ and $g (b, x) \in \psi_b (\alpha) \cap \psi_x (\beta)$.
\hfill$\square$

\bigskip

\section{Distributivity}

\medskip

\begin{Obs} Given a $\kappa$-complete, fine ideal $J$ on $P_\kappa (\lambda)$, the following are equivalent :
\begin{enumerate} [\rm (i)]
\item $J$ is $(\lambda^{< \kappa}, 2)$-distributive.
\item Given $A \in J^+$ and $F : A \times A \rightarrow \lambda^{< \kappa}$ with the property that $\vert \{F (a, b) : b \in A\} \vert < \kappa$ for all $a \in A$ , there is $B \in J^+ \cap P (A)$ and $h : B \rightarrow \lambda^{< \kappa}$ such that $\{ b \in B : F (a, b) \not= h (a)\} \in J$ for all $a \in B$.  
\item $J^+ \xrightarrow{J} (J^+)^n_\rho$ holds whenever $2 \leq n < \omega$ and $0 < \rho < \kappa$.
\item $J^+ \xrightarrow{J} (J^+)^3_2$ holds.
\end{enumerate}
\end{Obs} 

{\bf Proof.}  \hskip0,4cm  (i) $\rightarrow$ (ii) : Use Observation 5.4.

\medskip

\hskip0,2cm  (ii) $\rightarrow$ (iii) :  Assume that (ii) holds. It is readily seen that $J^+ \xrightarrow{J} (J^+)^2$ (and hence $J^+ \xrightarrow{J} (J^+, \kappa)^2$) holds. By Observation 5.18, it follows that $\kappa$ is weakly compact (and therefore inaccessible). Now to prove (iii), we proceed by induction on $n$. For $n = 2$, the assertion easily follows from (ii). Suppose now that the assertion has been verified for a certain $n$. Fix $A \in J^+$ and $F : A^{n + 1} \rightarrow \rho$, where $2 \leq \rho < \kappa$. For $b \in A$, let $T_b$ denote the collection of all functions $t$ from $(A \cap P (b))^{n - 1}$ to $\rho$. Notice that $\vert T_b \vert < \kappa$. Define $G : A \times A \rightarrow \bigcup_{b \in A} T_b$ as follows. Given $(b, c) \in A \times A$, let $G (b, c)$ be the element $t$ of $T_b$ defined by $t (a_1, \cdots, a_{n - 1}) = F (a_1, \cdots, a_{n - 1}, b, c)$. We may find $B \in J^+ \cap P (A)$, $h \in \prod_{b \in B} T_b$, and $X_b \in J$ for $b \in B$ such that $G (b, c) = h (b)$ whenever $b, c \in B$ and $c \notin X_b$. Define $H : \bigcup_{b \in B} ((A \cap P (b))^{n - 1} \times \{ b \}) \rightarrow \rho$ by $H (a_1, \cdots, a_{n - 1}, b) = h (b) (a_1, \cdots, a_{n - 1})$. There must $C \in J^+ \cap P (B)$, $i < 2$, and $Y_a \in J$ for $a \in C$ with 

\centerline{$\{b \in P_\kappa (\lambda) : a \setminus b \not= \emptyset\} \subseteq Y_a$}

such that $H (a_1, \cdots, a_{n - 1}, a_n) = i$ whenever $a_1, \cdots, a_n \in C$ and $a_{m + 1} \notin Y_{a_1} \cup \cdots \cup Y_{a_m}$ for $1 \leq m < n$. For $a \in C$, put $Z_a = X_a \cup Y_a$. Then clearly, $F (a_1, \cdots, a_{n + 1}) = G (a_n, a_{n + 1}) (a_1, \cdots, a_{n - 1}) = h (a_n) (a_1, \cdots, a_{n - 1}) = H (a_1, \cdots, a_{n - 1}, a_n) = i$ whenever $a_1, \cdots, a_{n + 1} \in C$ and $a_{m + 1} \notin Z_{a_1} \cup \cdots \cup Z_{a_m}$ for $1 \leq m \leq n$.

\medskip

\hskip0,2cm  (iii) $\rightarrow$ (iv) :  Trivial.

\medskip
\hskip0,2cm  (iv) $\rightarrow$ (i) :  Assume that (iv) holds. Then by Observation 5.18, $\kappa$ is inaccessible.  Now let $A \in J^+$, and $W_a \subseteq P_\kappa (\lambda)$ for $a \in P_\kappa (\lambda)$. Define $F : A \times A \times A \rightarrow 2$ by : $F (a, b, c) = 0$ if and only if $\{d \subseteq a : b \in W_d\} = \{d \subseteq a : c \in W_d\}$. There must be $B \in J^+ \cap P (A)$, $i < 2$, and $Z_a \in J$ for $a \in B$ such that $F (a, b, c) = i$ whenever $a, b, c \in B$, $b \notin Z_a$ and $c \notin Z_a \cup Z_b$. Given $a \in B$, we may find $C \in J^+ \cap P (B)$ such that $\{d \subseteq a : b \in W_d\} = \{d \subseteq a : c \in W_d\}$ whenever $d \subseteq a$ and $b, c \in C$. It easily follows that $i = 0$. Now fix $d \in P_\kappa (\lambda)$. Pick $a \in B$ with $d \subseteq a$, and $b \in B \setminus Z_a$. Then either $B \setminus (Z_a \cup Z_b) \subseteq W_d$, or $B \setminus (Z_a \cup Z_b) \subseteq P_\kappa (\lambda) \setminus W_d$.
%\hskip0,2cm  (ii) $\rightarrow$ (i) :  Suppose that (ii) holds, and let $A \in J^+$, and $W_a \subseteq P_\kappa (\lambda)$ for $a \in P_\kappa (\lambda)$. For $a \in P_\kappa (\lambda)$, set $W^{0}_a = W_a$ and $W^{1}_a = P_\kappa (\lambda) \setminus W_a$.  It is readily seen that $J^+ \xrightarrow{J} (J^+)^2$ (and hence $J^+ \xrightarrow{J} (J^+, \kappa)^2$) holds. By Observation, it follows that $\kappa$ is weakly compact (and therefore inaccessible).  Hence we may find $B \in J^+ \cap P (A)$ and $h_a : P (a) \rightarrow 2$ for $a \in B$ such that \centerline{$\{ b \in B : \exists c \subseteq a (b \notin W_c^{h_a (c)})\} \in J$} for all $a \in B$. Put $h = \bigcup_{a \in B} h_a$. Then, as is easily seen, $h$ is a func. tion from $P_\kappa (\lambda)$ to $2$. Moreover, $B \setminus W_c^{h (c)} \in J$ for any $c \in P_\kappa (\lambda)$.
\hfill$\square$

\medskip

Note the similarity with Fact 7.3 (ii) or Fact 7.12. One could indeed argue that $(\lambda^{< \kappa}, 2)$-distributivity of $J$ (respectively, mild $\lambda^{< \kappa}$-ineffability of $\kappa$) makes more sense (or is more natural) than $(\lambda, 2)$-distributivity (respectively, mild $\lambda$-ineffability). However the question of Abe mentioned in the introduction concerns $(\lambda, 2)$-distributivity (and not  $(\lambda^{< \kappa}, 2)$-distributivity) of $I_{\kappa, \lambda}$, so let us focus on $(\lambda, 2)$-distributivity.  

\medskip

\begin{Pro} Suppose that $J$ is a $(\kappa, 2)$-distributive, $\kappa$-normal, fine ideal on $P_\kappa (\lambda)$. Then $J \upharpoonright \kappa$ is $(\kappa, 2)$-distributive. 
\end{Pro}

{\bf Proof.}  Let $X \in (J \upharpoonright \kappa)^+$,  and $W_\alpha \subseteq \kappa$ for $\alpha < \kappa$. For $\alpha < \kappa$, set $T_\alpha = \{a \in \Omega_{\kappa, \lambda} : aÊ\cap \kappa \in W_\alpha\}$. We may find $h \in \prod_{\alpha < \kappa} \{T_\alpha, P_\kappa (\lambda) \setminus T_\alpha \}$, $B \in J^+ \cap P (\{a \in \Omega_{\kappa, \lambda} : a \cap \kappa \in X\})$, and $Z_\alpha \in J$ for $\alpha < \kappa$ such that $B \setminus Z_\alpha \subseteq h (\alpha)$ for all $\alpha < \kappa$. Define $k \in \prod_{\alpha < \kappa} \{W_\alpha, \kappa \setminus W_\alpha \}$ by : $k (\alpha) = W_\alpha$ if and only if $h (\alpha) = T_\alpha$. Put $S = \{a \in B : \forall \alpha \in a \cap \kappa (a \notin Z_\alpha)\}$ and $Y = \{a \cap \kappa : a \in S\}$. It is easy to see that $Y \in(J \upharpoonright \kappa)^+  \cap P (X)$. Furthermore $Y \setminus (\alpha + 1) \subseteq k (\alpha)$ for all $\alpha < \kappa$. 
 \hfill$\square$

\begin{Obs} Suppose that $\overline{\mathfrak{d}}_\kappa \leq \lambda = \cov (\lambda, \kappa^+, \kappa^+, \kappa)$. Then for any $D \in NS^\ast_{\kappa, \lambda}$ and any $X \in NS^+_\kappa \cap NCI_\kappa$, $I_{\kappa, \lambda} \vert \{a \in D : \sup (a \cap \kappa) \in X\}$ is not $(\kappa, 2)$-distributive. 
\end{Obs}

{\bf Proof.} Given $D \in NS^\ast_{\kappa, \lambda}$ and $X \in NS^+_\kappa \cap NCI_\kappa$, set $K = I_{\kappa, \lambda} \vert \{a \in D : \sup (a \cap \kappa) \in X\}$. Assume toward a contradiction that $K$ is $(\kappa, 2)$-distributive. By Fact 5.31, there is $C \in NS^\ast_{\kappa, \lambda}$ such that $I_{\kappa, \lambda} \vert C$ is $\kappa$-normal. Put $J = K \vert C$. Then clearly, $J$ is $\kappa$-normal and $(\kappa, 2)$-distributive. Hence by Observation 5.21 and Proposition 8.2, $J \upharpoonright \kappa$ is a normal, $(\kappa, 2)$-distributive, fine ideal on $\kappa$. It follows that $X \in J \upharpoonright \kappa$. Contradiction.
 \hfill$\square$

\begin{Pro} Suppose that $\overline{\mathfrak{d}}_\kappa \leq \lambda$ and $\cf (\lambda) \not= \kappa$. Then for any $D \in NS^\ast_{\kappa, \lambda}$ and any $X \in NS^+_\kappa \cap NCI_\kappa$, $I_{\kappa, \lambda} \vert \{a \in D : \sup (a \cap \kappa) \in X\}$ is not $(\lambda, 2)$-distributive. 
\end{Pro}

{\bf Proof.}  By Fact 7.4, Proposition 7.12 and Observation 8.3.
 \hfill$\square$

\bigskip

\section{The case $\cf (\lambda) = \kappa$}

\medskip

By Fact 4.4, if $\cf (\lambda) = \kappa$, then $\cov (\lambda, \kappa^+, \kappa^+, \kappa) > \lambda$. Thus a different approach is needed in case $\cf (\lambda) = \kappa$.

\medskip

\begin{Pro} Let $J$ be a $\kappa$-complete, fine ideal on $P_\kappa (\lambda)$ such that $cof (J) < \max \{\mathfrak{d}_\kappa, u (\kappa^+, \lambda)\}$. Then the following hold : 
\begin{enumerate} [\rm (i)]
\item Let $A_\alpha \in J^+$ for $\alpha < \kappa$ be given such that $A_\alpha \subseteq A_\beta$ whenever $\beta < \alpha < \kappa$. Then there is $C \in J^+$ such that $C \setminus A_\alpha \in I_{\kappa, \lambda}$ for all $\alpha < \kappa$. 
\item Suppose that $\kappa$ is weakly compact. Then given $W_\alpha \subseteq P_\kappa (\lambda)$ for $\alpha < \kappa$, there is $C \in J^+$ and $h \in \prod_{\alpha < \kappa} \{ W_\alpha, P_\kappa (\lambda) \setminus W_\alpha\}$ such that $C \setminus h (\alpha) \in I_{\kappa, \lambda}$ for all $\alpha < \kappa$. 
 \end{enumerate}
\end{Pro} 

{\bf Proof.}  (i) : The proof is a straightforward modification of that of Proposition 2.7 in \cite{MP}.

(ii) :  Assume that $\kappa$ is weakly compact, and let $W_\alpha \subseteq P_\kappa (\lambda)$ for $\alpha < \kappa$. There must be $h \in \prod_{\alpha < \kappa} \{ W_\alpha, P_\kappa (\lambda) \setminus W_\alpha\}$ such that $A_\alpha \in J^+$ for all $\alpha < \kappa$, where $A_\alpha = \bigcap_{\beta \leq \alpha} h (\beta)$. By (i) we may find $C \in J^+$ such that $C \setminus A_\alpha \in I_{\kappa, \lambda}$ for all $\alpha < \kappa$. Then clearly, $C \setminus h (\alpha) \in I_{\kappa, \lambda}$ for every $\alpha < \kappa$.
\hfill$\square$
 
\begin{Cor} Suppose that $\kappa$ is weakly compact, and let $J$ be a $\kappa$-complete, fine ideal on $P_\kappa (\lambda)$ such that $cof (J) < \max \{\mathfrak{d}_\kappa, u (\kappa^+, \lambda)\}$. Then the following hold : 
\begin{enumerate} [\rm (i)]
\item $J$ is $(\kappa, 2)$-distributive. 
\item $J^+ \xrightarrow[\kappa]{I_{\kappa, \lambda}} (J^+)^2_\rho$ holds for any nonzero cardinal $\rho < \kappa$.
 \end{enumerate}
\end{Cor} 

{\bf Proof.} Use Observation 5.4. 
\hfill$\square$

\begin{fact}  {\rm(\cite{square}, \cite{MS})} Suppose that $\kappa$ is weakly inaccessible, and let $J$ be a $\kappa$-complete, fine ideal on $P_\kappa (\lambda)$ such that $cof (J) < \max \{\mathfrak{d}_\kappa, u (\kappa^+, \lambda)\}$. Then $J^+ \xrightarrow[\kappa]{I_{\kappa, \lambda}} [J^+]^2_{\kappa^+}$ holds.
\end{fact} 

\medskip

To obtain a negative partition relation, we will go one cardinal up and work with partitions of $\kappa^+ \times P_\kappa (\lambda)$. 

\medskip

\begin{Def} Given $2 \leq n < \omega$, two collections $X$ and $Y$ of subsets of $P_\kappa (\lambda)$, an ideal $J$ on $P_\kappa (\lambda)$, and a cardinal $\rho$, $X \xrightarrow[\kappa^+]{Y} (J^+)^2_\rho$ means that for any $F : \kappa^+ \times P_\kappa (\lambda) \rightarrow \rho$ and any $A \in X$, there is $i < \rho$ and $B \in J^+ \cap P (A)$ such that $\{ b \in B : F (\sup (a \cap \kappa^+), b) \not= i\} \in Y$ for all $a \in B$.  
\end{Def}

\begin{fact} {\rm (\cite{square})} Suppose that the following hold :
\begin{itemize}  
\item $\kappa$ is weakly inaccessible.
\item $\cf (\lambda) = \kappa$.
\item $2^\kappa \leq \lambda$.
\item $u (\kappa^+, \tau) \leq \lambda$ for any cardinal $\tau$ with $\kappa < \tau < \lambda$.
\end{itemize}
Then $\{C \cap S^\lambda_{\kappa, \lambda}\} \xrightarrow[\kappa^+]{I_{\kappa, \lambda}} [I_{\kappa, \lambda}^+]^2_\lambda$ fails for some $C \in NS^\ast_{\kappa, \lambda}$.
\end{fact}

 \begin{Pro} Suppose that $\cf (\lambda) = \kappa$ and $2^\kappa \leq \lambda$. Let $D \in NS^\ast_{\kappa, \lambda}$ and $J = I_{\kappa, \lambda} \vert (D \cap S^\lambda_{\kappa, \lambda})$. Then the following hold : 
 \begin{enumerate} [\rm (i)]
 \item  $J^+ \xrightarrow{J} (J^+)^2_\omega$ does not hold.
 \item $J$ is not $(\lambda, 2)$-distributive.
 \end{enumerate}
\end{Pro} 

{\bf Proof.}  (i) : By Observations 5.18 and 6.4 and Fact 9.5.

(ii) : Suppose otherwise. Then $\lambda^{< \kappa} = \lambda$ by Fact 7.10, since $\kappa$ is mildly $\lambda$-ineffable by Fact 7.4. Hence $J^+ \xrightarrow{J} (J^+)^2_\omega$ holds by Observation 8.1, which contradicts (i).
\hfill$\square$

  \bigskip
\noindent Universit\'e de Caen - CNRS \\
Laboratoire de Math\'ematiques \\
BP 5186 \\
14032 Caen Cedex\\
France\\
Email :  pierre.matet@unicaen.fr\\

\end{document}